\RequirePackage{fix-cm}
\documentclass[smallcondensed]{svjour3}       % onecolumn (second format)
\smartqed  % flush right qed marks, e.g. at end of proof
\usepackage{mathptmx}      % use Times fonts if available on your TeX system
\usepackage{amsmath}

\usepackage{amsthm}
\usepackage{amsfonts}
\usepackage{fancybox}
\usepackage{amssymb}
\usepackage{epic,eepic}
\usepackage{latexsym,bm}
\usepackage{float}
\usepackage{color}
\usepackage{enumerate} 
\usepackage{extarrows}
\usepackage{graphicx,epstopdf,epsfig,epsf}
\graphicspath{{figs/}}
\usepackage[colorlinks,linkcolor=blue,citecolor=blue,urlcolor=blue]{hyperref}
\usepackage{algorithm}
\usepackage{algorithmic}
\usepackage{cite}

\newcommand{\R}{\mathbb{R}}
\newcommand{\E}{\mathbb{E}}

\newcommand{\I}{\mathbb{I}}
\newtheorem{assumption}{Assumption}

\newtheorem{lem}{Lemma}
\newtheorem{thm}{Theorem}

\allowdisplaybreaks[4]
\numberwithin{equation}{section}
\numberwithin{algorithm}{section}
\newtheorem{rem}{Remark}

\newtheorem{defi}{Definition}

\numberwithin{rem}{section}
\numberwithin{thm}{section}
\numberwithin{defi}{section}
\numberwithin{ass}{section}
\numberwithin{lem}{section}
\usepackage[misc]{ifsym}
\usepackage[all]{hypcap}

\begin{document}
	
	\title{Stochastic Trust-Region Methods with Trust-Region Radius Depending on Probabilistic Models 
	}

	
	\author{Xiaoyu Wang     \and
		Ya-xiang Yuan %etc.
	}
	
	
	\institute{ Xiaoyu Wang(\Letter) \at
		Institute of Computational
		Mathematics and Scientific/Engineering Computing, Academy of Mathematics and Systems Science, Chinese Academy of Sciences, Zhong Guan Cun Donglu 55, Beijing 100190, China. \& University of Chinese Academy of Sciences, Beijing 100049, China. \\
		\email{wxy@lsec.cc.ac.cn}           %  \\
		\and
		Ya-xiang Yuan \at
		State Key Laboratory of Scientific/Engineering Computing, 
		Institute of Computational
		Mathematics and Scientific/Engineering Computing, Academy of Mathematics and Systems Science, Chinese Academy of Sciences, Beijing 100190, China. \\
		\email{yyx@lsec.cc.ac.cn} 
	}
	
	\date{Received: date / Accepted: date}

	\maketitle
	
	\begin{abstract}
		We present a stochastic trust-region model-based framework in which its radius is related to the probabilistic models.  Especially, we propose a specific algorithm, termed STRME, in which the trust-region radius depends linearly on the latest model gradient. The complexity of STRME method in non-convex, convex and strongly convex settings has all been analyzed, which matches the existing algorithms based on probabilistic properties. In addition, several numerical experiments are carried out to reveal the benefits of the proposed methods compared to the existing stochastic trust-region methods and other relevant stochastic gradient methods. 
	\end{abstract}
	
	\section{Introduction}
	\label{intro}
	
	In this paper, we are concerned with the following unconstrained optimization problem
	\begin{equation}
	\label{Fun1}
	\min_{x \in \R^d} \, f(x),
	\end{equation}
	where the objective function $f$ is assumed to be smooth and bounded from below.
	But we only have access to the value of $f$ and its derivative information with some noise. In recent years, the expected risk minimization (ERM) problem, which is fundamental in the field of machine learning and statistic, has become the focus of many researchers. The ERM problems can be formulated as follows:
	\begin{equation}
	\label{Fun2}
	\min_{x\in \R^d}\, f(x) = \E [ f(x, \xi)],
	\end{equation}
	where $\E[\cdot]$ denotes the expectation taken with respect to the random variable $\xi \in \R^d$.
However, because the probability distribution of $ \xi$ is unknown in advance, solving (\ref{Fun2}) is intractable directly. Usually only noisy information about the gradient of $f$ is available.  The empirical risk problem with a fixed amount of data (possibly very large) or the on-line setting problem where the data is flowing in sequentially, which involves an estimate of problem ($\ref{Fun2}$), is more often considered in practice. Through the whole paper, we mainly consider stochastic optimization methods to  solve such kind of problems.
	
	The classic stochastic optimization method is stochastic gradient descent (SGD) method, which dates back to the work by Robbins and Monro \cite{SGD-1951}. The method is prominent and adorable in large-scale machine learning due to simpleness and low-cost computing. However, because of the variance introduced by random sampling, the sequence of learning rate (step-size) progressively diminish both in theoretical analysis and practical implementation, which leads to slow convergence. Thus finding an appropriate learning rate is critical for the performance of SGD method, but it is not easy in practice. To deal with aforementioned issues, various adaptive gradient algorithms have emerged, for instance AdaGrad\cite{AdaGrad}, RMSProp\cite{RMSProp}, Adam\cite{Adam}, which are very popular in deep learning. 
	Besides, variance reduction (VR) methods, to improve the performance of SGD method, are proposed, such as SVRG\cite{SVRG}, SAGA\cite{SAGA} and SARAH \cite{SARAH}. Especially, they have achieved linear convergence rate when solving strongly convex problems, which is a stronger result than that of SGD method. Furthermore, these methods have also been extended to solve non-convex problems such as the deep neural networks, and achieve good performance \cite{SAGA-nonconvex, SVRG-nonconvex, SARAH_nonconvex}. The VR technique is applicable to the problem with a large but fixed sample set, for which the full gradient has to be calculated as a compromise to achieve the significant variance reduction. Hence they are not easy to fit on the on-line setting like SGD method and adaptive gradient methods.
	
	 Besides, many second order methods are proposed, which are known to perform better than the first-order methods on various highly nonlinear and ill-conditioned problems\cite{oLBFGS, SQN, SdLBFGS, NIM, IQN}. More recently, cubic regularization methods as a class of Newton-type variants have attracted a lot of interest\cite{Xu_Peng,CR_Kohler,CR_Ghadimi,CR_Jordan,CR_VR_Lan,CR_VR_Zhou}. Especially, the VR technique is applied to improve the performance of cubic regularization methods \cite{CR_VR_Lan,CR_VR_Zhou}.
 
	 Sample averaging is a natural and well-known technique to reduce the variance of gradient or the noise\cite{Sample_average}. Generally speaking, to guarantee the accuracy of the function and gradient estimators, the sample size has to increase when the algorithm goes to optimality. And the training sample is not only regarded as a fixed and finite set, that is to say that the sample averaging technique can be employed to the on-line setting. %\red{In our numerical experiments, we will use the sample averaging technique to estimate the function value and gradient of $f$.} 

	Recently, with the success of deep neural networks, the development and analysis of methods for non-convex problems have attracted tremendous attention. As we know, traditional trust-region methods is a class of well-established and effective methods in nonlinear optimization\cite{Powell_TR_a, Yuan_TR_review}. For details interested readers can refer to the review by Yuan\cite{Yuan_TR_review}. With such a framework, we can utilize second order information when building the trust-region subproblem. Besides, due to the boundedness of the trust-region, the Hessian approximation matrix is not required to be positive definite. An advantage of trust-region methods is that they can be applied to non-convex and ill-conditioned problems. Although there are various effective methods as mentioned before, trust-region algorithms deserve more attention in stochastic optimization.  
	
	Actually, the traditional trust-region framework has already been considered to solve machine learning problems \cite{TR_logistic, TR_Newton_linear_classification, SFN_saddle, Two_stage_TR, Xu_Peng}. Dauphin et al. \cite{SFN_saddle} proposed a saddle free Newton (SFN) which exploits the exact Hessian information to escape saddle points. However its computation is high cost for large-scale and high-dimension problems. A two-stage subspace trust-region approach \cite{Two_stage_TR} was proposed to train deep neural networks, in which the local second-order model is conducted based on the partial information computed from a subset of the data. But the approach lacks theoretical guarantees. \cite{TR_Newton_linear_classification, TR_logistic} are designed to solve a specific class of machine learning problems, which need to utilize the accurate derivative information to construct good models, and the accurate function values to obtain good estimators, of which the computational costs is too expensive to afford for general large-scale machine learning problems. \cite{Xu_Peng} incorporates inexact Hessian information into the trust-region framework but the exact gradient and function values are required to be computed per iteration. 
	
	In \cite{Bandeira_TR_2014, Gratton_TR}, the authors construct the inexact models to satisfy some first-order accurate conditions with sufficiently high probability when building the trust-region subproblem. Cartis and Scheinberg \cite{Cartis_Linesearch} analyzes the complexity of line search and cubic regularization algorithms, which is based on the random models with certain probability but their function estimators are accurate, in non-convex, convex and strongly convex settings. These mainly focus on derivative free optimization (DFO) problems.

	A stochastic trust-region algorithm named STORM for stochastic optimization setting has been introduced in \cite{STORM}. Not only the model is conducted to satisfy some first-order accuracy requirements, but also the function values both at current iterate and next potential iterate are estimated with some probability, instead of the exact function values. And the liminf-type and lim-type first-order convergence has been analyzed. In addition, Blanchet et al. \cite{STORM_nonconvex} has bounded the expected convergence complexity of STORM for non-convex problems. More recently, Paquette and Scheinberg \cite{S_line_search} analyzes the complexity of a stochastic line search algorithm, of which the gradient and function estimators are both randomly sampled with some probability.

  Previous works have established the convergence and complexity properties for such a trust-region framework. In this paper, we are particularly interested in introducing a new trust-region radius formula which depends on the latest probabilistic model, to improve the practical performance of such trust-region methods. Besides, we present an algorithm termed STRME in which the trust-region radius depends linearly on the model gradient just updated, following a piece of work initially proposed for deterministic optimization\cite{Fan_Yuan}. The idea is meaningful and attractive. Because the trust region can be tailored according to the newly generated model, not just the success of the trial steps as that of STORM. Note that the trust-region radius of STORM at iteration $k$ is completely determined by the past iterations $0$ to $k-1$. However, the framework we proposed is no longer that case. That is to say the trust-region radius in STRME, which related to the current model, is not measurable with respect to all the information generated from the past iterations $0$ to $k-1$.
  	
  	It will bring new challenges both in theoretical analysis and numerical experiments. The trust-region radius in STRME depends on the newly updated model makes it more complicated to analyze the complexity of STRME, compared to that of STORM. In \cite{S_line_search}, Paquette and Scheinberg tackles this issue, of which the quantity $\Delta_k$ is also related to the currently updated gradient. Note that the quantity $\Delta_k$ in \cite{S_line_search} is not really the trust-region radius of trust-region algorithms. To obtain the complexity results of STRME, we analyze the parameter $\Lambda_k$ (called relative trust-region radius), which depends on the past iterative information, instead of the trust-region radius itself. Our approach avoids the difficulty of analyzing the trust-region radius directly. 
	The convergence analysis for such trust-region algorithms relies on the requirements that these quantities such as models and function estimators are sufficiently accurate with sufficiently high probability. And the accuracy of these quantities are controlled by the trust-region radius. However the trust-region radius is unknown before the model is updated. We will elaborate this at the beginning of Section \ref{sec:4} where we present numerical experiments.
	
	These changes indeed bring forth some advantages. The choice of trust-region radius makes STRME algorithm scale invariant on problems. Besides, the trust-region radius which depends on the model gradient can capture more new information. We have to say that it is advisable to adjust the trust-region radius according to the latest probabilistic model. Our numerical experiments illustrate this viewpoint. We have tested on regularized logistic regression and a simple deep neural network on real datasets. The numerical experiments show that the trust-region radius of the proposed STRME method reduces asymptotically. And the oscillation in the trust-region radius is less severe in contrast to that of STORM method. In addition, we observe that the proposed algorithm can get more successful iterates after a long time training.

	The major contributions of this paper are summarized as follows:
	\begin{itemize}
		\item[(1)] The trust-region radius can be defined using the probabilistic model, in particular its gradient.
		\item[(2)] We propose a specific algorithm termed STRME, in which the trust-region radius $\delta_k$ is linearly dependent on the norm of the model gradient. The complexity of STRME in non-convex, convex and strongly convex cases are analyzed, respectively. The expected number of iterations of STRME algorithm for non-convex problem is $O(\epsilon^{-2})$ by reaching $\left\|\nabla f(x) \right\|\leq \epsilon $, which is similar to the result in \cite{STORM_nonconvex}.		
		In addition, the expected convergence rates for general convex and strongly convex problems, which are $O(1/\epsilon)$ and $O(\log(1/\epsilon))$ for reaching $f(x) - f^{\ast} \leq \epsilon$, respectively. 
		\item[(3)] In numerical experiments, sample averaging technique is utilized to construct probabilistic models. Besides, we adopt the dogleg method to solve the trust-region subproblem for the regularized logistic regression problem. In the same way, the limited memory symmetric rank one (L-SR1) is employed to approximate the Hessian matrix, and then incorporate them into STRME algorithm to train a deep neural network problem. The results indicate that the proposed algorithm compares favorably to other stochastic optimization algorithms.

	\end{itemize}
	
	The outline of this paper is as follows. In Section 2 we give some definitions about the probabilistic models and estimates, and present a generic analysis framework based on the random models and estimates; In Section 3 we propose a specific algorithm named STRME and analyze the complexity in non-convex, convex and strongly convex cases; In Section 4 we report some numerical results on regularized logistic regression problem and a simple deep neural network to show the efficiency of STRME in different settings; In the end, we draw some conclusions in Section 5.
	
	\paragraph{	\textbf{Notations.}}
	%The class of all functions with the form of (1.2) is denoted by $\mathcal{F}_n$.
	Throughout this paper, we use $x^{\ast}$ to denote the global minimizer, $f^{\ast}=f(x^{\ast})$. Let $\left\|\cdot \right\| $ denote the Euclidean norm, i.e. $\left\|\cdot\right\|_2 $, unless otherwise specified. Let $B(x, \Delta)$ denote the ball of the radius $\Delta$ around $x$.  Let $\I\left\lbrace A\right\rbrace $ denote the indicator function of the event $A$, that is: if $A$ occurs, $\I\left\lbrace A\right\rbrace =1 $; else, $\I\left\lbrace A\right\rbrace = 0.$ A function $f \in \mathcal{C}^{1}(\R^d)$, if the first derivation of $f$ exists and continuous. A function $f$ is $L$-smooth, if there is a constant  $L > 0$ such that $
	\left\|\nabla f(x) -\nabla f(y) \right\| \leq L \left\|x-y \right\|, \quad \forall \,x, y \in \R^d.$
	\section{A generic analysis framework based on random models and estimates}\label{sec:2}
	Let us first introduce a generic stochastic trust-region framework. The analysis for the framework can particularize to the specific algorithm for example the algorithm STRME  proposed in Section \ref{sec:3}, and the objective function, whether it is non-convex or convex, provided that the assumptions are satisfied.
	
	%%%%%%%%%%%%%%%%%%%%这里算法1和2，算法1报告STORM，然后算法2作为1的特例
	\begin{algorithm}[H]
		\begin{algorithmic}[1]
			\STATE \textbf{Initialization}: Given an initial point $x_0$, $\gamma > 1$, $\eta_1\in(0,1)$, $\eta_2 > 0$ $\mu_0 \in (0, \mu_{\max})$ with $\mu_{\max} > 0$;  Set $k = 0$	
			\STATE Construct a model(possibly random) $m_k(x_k + d)$ to approximate $f(x)$ at $x_k$ with $d = x - x_k$
			\STATE Compute $\delta_{k} = \delta(m_k, \mu_k) $
			\STATE Compute a trial step $d_k = \arg\min_{\left\| d\right\| \leq \delta_k} m_k(x_k + d)$ to such that $d_k$ satisfies a sufficient reduction condition
			\STATE Obtain estimates $f_k^0$ and $f_k^d$ of $f(x_k)$ and $f(x_k + d_k)$
			\STATE Compute $\rho_k = \frac{f_k^0 - f_k^d}{m_k(x_k) - m_k(x_k + d_k)}$
			\IF{$ \rho_k \geq \eta_1$ and $\left\|g_k \right\|\geq \eta_2\delta_k $}
			\STATE$x_{k+1} = x_k + d_k$, $\mu_{k+1} = \min(\gamma\mu_k, \mu_{\max})$
			\ELSE 
			\STATE $x_{k+1} = x_k$, $\mu_{k+1} = \mu_k/\gamma$	
			\ENDIF	
			\STATE Set $k := k+1$, and go to step 2
			
		\end{algorithmic}
		\caption{A Stochastic Trust Region Framework  }\label{alg:1} 
	\end{algorithm}

	The proposed Algorithm \ref{alg:1} covers the framework of STORM \cite{STORM}. We can see that if $\delta_k = \mu_k$, the above algorithm will reduce to STORM algorithm. The main difference lies in the trust-region radius in Algorithm \ref{alg:1} which depends on the current model $m_k$. We directly update the parameter $\mu_k$ (called relative trust-region radius), not the trust-region radius $\delta_k$ itself, which avoids the troubles by the randomness rise to the currently random model. Of course, except the current random model, there may be other factors, {for instance the iterative models of previous steps, i.e. $\left\lbrace m_t\right\rbrace_{t\leq k} $, and the difference of previous iterates $\left\lbrace x_{t}-x_{t-1}\right\rbrace_{t\leq k} $ and so on. The framework we present here does not involve the specific forms of the random models and the sufficient reduction condition. We will discuss them in the next part.
	
	Note that Algorithm \ref{alg:1} generates a random process. Obviously, the randomness of the algorithm comes from the randomness of the models and estimates we have constructed per iteration. At a deep level, it is determined by the inexact or random information obtained from the problem we are trying to solve. To formalize the random process, we introduce some notations to describe the quantities of them and their realizations. Let $M_k$ denote the random model in $k$-th iteration, while $m_k = M_k(\omega)$ for its realization, where $\omega$ is a random variable. We know that the randomness of the models gives rise to the randomness of the iterates, relative trust-region radius and trial step produced by algorithm \ref{alg:1}. These random variables are denoted by $X_k$, $\Lambda_k$ and $D_k$, respectively, while let  $x_k = X_k(\omega)$, $\mu_k = \Lambda_k(\omega)$, and $d_k = D_k(\omega)$ to denote their realizations. Similarly, we use $\left\lbrace F_k^0, F_k^d \right\rbrace $ to denote the random estimates of $f(X_k)$ and $f(X_k + D_k)$, while their realizations are denoted by $f_k^0 = F_k^0(\omega)$ and $f_k^d = F_k^d(\omega)$. We will utilize those notations to analyze the random process later in this section, which is under some assumptions that model $M_k$ and estimates $F_k^0, F_k^d$ are sufficiently accurate with some probability conditioned on the past. In order to formalize all the randomized information before $k$-th iteration, let $\mathcal{F}_{k-1}^{M\cdot F}$ denote the $\sigma$-algebra generated by $\left\lbrace M_0, \cdots, M_{k-1}\right\rbrace $ and $\left\lbrace F_0^0, F_0^d, \cdots, F_{k-1}^0, F_{k-1}^{d}\right\rbrace $.  After the current model $m_k$ is constructed, let $\mathcal{F}_{k-1/2}^{M\cdot F}$ denote the $\sigma$-algebra generated by $\left\lbrace M_0, \cdots, M_{k}\right\rbrace $, and $\left\lbrace F_0^0, F_0^d, \cdots, F_{k-1}^0, F_{k-1}^{d}\right\rbrace $.	
	
	Next, we will introduce some definitions to precise our requirements on the probabilistic models and estimates.
	
	%%%%%%%%%%%%%%%%%%%%%%%%%%%%%%%%%%%%%%%%%%%%%%%%%%%%%%%%%
	%%%%%%%%%%%%%%%%%%%%%%%%%%%%%%%%%%%%%%%%%%%%%%%%%%%%%%%%%		
	\subsection{Probabilistic models and estimates}	
	\label{sec:2:1}
	
	First we recall the measure for accuracy of deterministic models, which is introduced in \cite{Random_model_1,Random_model_2}.

	\begin{defi}\label{defi_det_model}
		We say a model $m_k$ is $\kappa$-fully linear model of f on $B(x_k,\delta_k)$, for $\kappa = (\kappa_{ef},\kappa_{eg})$,  if \,  $\forall y\in B(x_k,\delta_k)$, 
		\vspace{-0.5em}
		\begin{equation}
		\begin{split}
		\left\|\nabla f(y) - \nabla m_k(y) \right\|  & \leq  \quad \kappa_{eg}\delta_k, \,\,\,  \text{and} \\
		\left|f(y) - m_k(y) \right|   & \leq \quad  \kappa_{ef} \delta_k^2. 	 		
		\end{split}
		\end{equation} 
	\end{defi}
	The extending concept of the above definition is probabilistically fully-linear model which is described in \cite{STORM}.
	\begin{defi}\label{Random_model}
		A sequence of random model $M_k $ is said to be $\alpha$-probabilistically $\kappa$-fully linear
		with respect to the corresponding sequence  $\left\lbrace X_k, \Lambda_k\right\rbrace $, if the events
		\vspace{-0.5em}
		$$ I_k = \I{\left\lbrace M_k \text{ is a } \kappa\text{-fully linear model of} \,f \,\text{on}\,\, B(x_k, \delta(\mu_k))\right\rbrace} $$ 
		satisfy the condition: 
		\vspace{-0.5em}
		\begin{equation}
		P(I_k = 1 | \mathcal{F}_{k-1}^{M\cdot F}) \geq \alpha.
		\end{equation} 
	\end{defi}
	
	The above definition states that the model $M_k$ is a locally good approximation of the first-order Taylor expansion of the objective function with probability at least $\alpha$, conditioned on $\mathcal{F}_{k-1}^{M\cdot F}$. However, there is still some possibility such that the model is inaccurate, even very bad. To guarantee the quality of the trial step, we  hope the random model $M_k$ closer to the first-order Taylor expansion. However, the corresponding computation cost will increase. Thus there is a trade-off between the accuracy of the model and the computation cost.
	
	Taking aside of the accurate model, the estimates of $f(x_k)$ and $f(x_k + d_k) $ are also required to be sufficiently accurate. The deterministic version of accurate estimates is formally stated as follows. 
		
	\begin{defi}
		The estimates $f_k^0$ and $f_k^d$ are  $\epsilon_F$-accurate estimates of $f(x_k) $ and $f(x_k + d_k)$, if
		\vspace{-0.5em}
		\begin{equation}
		\begin{split}
		\left| f_k^0 - f(x_k) \right|  & \leq  \quad \epsilon_F
		\delta_k^2, \,\,\,  \text{and}   \\
		\left|f_k^d - f(x_k + d_k) \right|   & \leq \quad  \epsilon_F \delta_k^2. 	 		
		\end{split}
		\end{equation} 
	\end{defi}
	The definition of probabilistically accurate estimates is shown as follows which is a modified version of that in \cite{Random_estimate}. 
	\begin{defi}\label{Random_estimate}
		A sequence of random estimates $\left\lbrace F_k^0, F_k^d \right\rbrace  $ is said to be $\beta$-probabilistically $\epsilon_F$-accurate with respect to the corresponding sequence $\left\lbrace X_k, \Lambda_k, D_k \right\rbrace $, if the events
		\vspace{-0.5em}
		$$ J_k = \I\left\lbrace F_k^0 \,\text{and} \, F_k^d\, \text{are}\, \,\epsilon_F\text{-accurate estimates of} \,f(x_k)\, \text{and} \,f(x_k + d_k), \text{respectively} \right\rbrace  $$ 
		satisfy the condition: 
		\vspace{-0.5em}
		\begin{equation}
		P(J_k = 1 | \mathcal{F}_{k-1/2}^{M\cdot F}) \geq \beta,
		\end{equation}
		where $\epsilon_{F}$ is a fixed constant.
	\end{defi}	
	
	Using Definitions \ref{Random_model} and \ref{Random_estimate}, we assume that model $M_k$  and estimates $\left\lbrace F_k^0,F_k^d\right\rbrace $ satisfy the following assumption in our analysis.	
	
	\begin{assumption}\label{section2_assump1}
		The followings hold for the quantities used in Algorithm \ref{alg:1}
		\begin{enumerate}[(i)]
			\item \label{assump1(1)}There exist $\kappa_{ef}$,  $\kappa_{eg} > 0$ such that the sequence of random models $M_k$ is $\alpha$-probabilistically $(\kappa_{ef}, \kappa_{eg})$-fully linear, for a sufficiently large $\alpha \in (0,1)$. 
			\item\label{assump1(2)} There exist $\epsilon_F > 0$ such that
			the sequence of estimates $\left\lbrace F_k^0, F_k^d\right\rbrace $  is	 
			$\beta$-probabilistically $\epsilon_{F}$-accurate, for a sufficiently large $\beta \in (0,1)$.
			\item\label{assump1(3)} The sequence of estimates  $\left\lbrace F_k^0, F_k^d\right\rbrace $ generated by Algorithm \ref{alg:1} satisfies the following condition that 
			\begin{equation}\label{AS1_3}
			\begin{split}
			\E[ \left|F_k^0 - f(X_k) \right|| \mathcal{F}_{k-1/2}^{M,F} ] \leq \kappa_{f} \Delta_k^2,\,\, \text{and} \\ 
			\E[ \left|F_k^d - f(X_k + D_k) \right|| \mathcal{F}_{k-1/2}^{M,F} ] \leq \kappa_{f} \Delta_k^2,
			\end{split}
			\end{equation}
			with $\kappa_f > 0$.
		\end{enumerate}		
	\end{assumption}
	
	\begin{rem}
		The above assumption is significant in the following convergence rates analysis. Compared to the Assumption 3.1 in \cite{STORM_nonconvex}, we know that  Assumption \ref{section2_assump1}(\ref{assump1(3)}) is additional. Nevertheless it is essential, which states that estimates $\left\lbrace F_k^0, F_k^d\right\rbrace $ can not be too worse in expectation in contrast to the true function values $\left\lbrace f(X_k), f(X_{k} + D_k)\right\rbrace $  when estimates $\left\lbrace F_k^0, F_k^d\right\rbrace $ are not sufficiently accurate. One may doubt that this additional condition is a little stronger than that of stochastic line search \cite{S_line_search}. One would agree that (\ref{AS1_3}) given above and (2.3) of \cite{S_line_search} are close in spirit. However, it does not possible to obtain one from the another. Actually, the term $\Lambda_{k}\left\|\nabla f(X_k) \right\|^2 $ in (2.3) of \cite{S_line_search} is hard to be calculated directly in the construction of estimates $\left\lbrace F_k^0, F_k^d \right\rbrace $. From later analysis of Theorem \ref{thm_Phi_decrease}, we can see that if Assumption \ref{section2_assump1}(\ref{assump1(3)}) is relaxed as Assumption 2.4 in \cite{S_line_search}, the bound for $\E[\Phi_{k+1} - \Phi_{k}]$ will be worse than the current results. Moreover the estimates for $\nu$ and other parameters will be more complicated.  Note that Assumption \ref{section2_assump1}(\ref{assump1(3)}) can easily be satisfied in the practical implementation as long as Assumption\ref{section2_assump1}(\ref{assump1(2)}) holds. We will explain this later in our numerical experiments.
	\end{rem}	
	
	To make the analysis simple and easy to understand, we use the following statements.
	\begin{itemize}
		\item If $I_k = 1$, we say that the model is {\bf{true}}; otherwise, we say that the model is {\bf{false}}.
		\item If $J_k = 1$, we say that the estimates are {\bf{tight}}; otherwise, we say that the estimates are {\bf{loose}}.
		\item If an iteration $k$ is accepted, we say that the iteration is {\bf{successful}}; otherwise, we say that the iteration is {\bf{failed}}. 
	\end{itemize}		
	
	In the end of this subsection, we would like to introduce a definition of convergence criterion for analysis, which is named $\epsilon$-solution. When $f$ is unknown to be convex, we say that $X_k$ is an $\epsilon$-solution if $ \left\| \nabla f(X_k) \right\| \leq \epsilon$.
	However, when $f$ is convex or strongly convex, we say that $X_k$ is an $\epsilon$-solution if $ f(X_k) - f^{\ast}  \leq \epsilon$.	
	
	\begin{rem}
		There are three cases of non-convex, convex and strongly convex to be discussed in this article. Due to the intractability of the general non-convex problem, it is unreasonable to use the same criterion as the convex problem. Thus the definition of $\epsilon$-solution for non-convex case is different from that for convex case.
	\end{rem}

	%%%%%%%%%%%%%%%%%%%%%%%%%%%%%%%%%%%%%%%%%%%%%%%%%%%%%%%%%
	%%%%%%%%%%%%%%%%%%%%%%%%%%%%%%%%%%%%%%%%%%%%%%%%%%%%%%%%%
	\subsection{Analysis of the stochastic process}	
	\label{sec:2:2}
	In this part, we aim to estimate the upper bound of an expected stopping time by observing the behavior of stochastic process generated by Algorithm \ref{alg:1}. The results can be applied to analyze the convergence rates of the proposed algorithm in different settings.
	
	We first give some basic definitions before theoretical analysis. 
	\begin{defi}
		Let $X = \left\lbrace X_k, k \geq 0\right\rbrace $ be a stochastic process. We say $T$ is a stopping time with respect to $X$ if for each $k > 0$, the event $\left\lbrace T = k \right\rbrace $ is completely determined by the total information up to time $k$, that is $\left\lbrace X_0, X_1, \cdots,X_k\right\rbrace $.
	\end{defi}
	Here we give a random variable $T_{\epsilon}$, which is the total number of iterations until an $\epsilon$-solution is achieved. 
	\begin{rem}
		$T_{\epsilon}$ is a special stopping time and dependent on randomness of the proposed algorithm and the $\epsilon$-solution we have defined.
	\end{rem}
	
	Next, we consider a stochastic process $\left\lbrace \Lambda_k, \Phi_k\right\rbrace $ such that $\Lambda_k \in \left[ 0,\infty \right) $ and $\Phi_k \in \left[0, \infty \right) $ for all $k > 0$. Let us give the definition of a special random event $W_k$ as follows
	\begin{equation}\label{W_k}
	P(W_k = 1 | \mathcal{F}_{k-1}^{M,F}) = p, \,\, P(W_k = -1 | \mathcal{F}_{k-1}^{M,F}) = 1-p,
	\end{equation}
	where $p\in [0,1]$.
	We now assume that $\Lambda_k $ and $\Phi_k$ satisfy the assumption mentioned below, for all $ k < T_{\epsilon}$.
	
	\begin{assumption}\label{section2_assump2}
		
		\begin{enumerate}[(i)]
			\item[]		
			\item\label{assump2(1)} There exist constants $\Phi_{\max} > 0 $ and $\mu_{\max} > 0$ such that $\Phi_{k} \leq \Phi_{\max}$ and $\Lambda_k \leq \mu_{\max}$, respectively.
			
			\item\label{assump2(2)} There exists a constant $\hat{\Lambda} > 0$ such that for all $k \leq T_{\epsilon}$, the following properties hold  
			\begin{equation}
			\Lambda_{k+1} \geq \min(\Lambda_ke^{\lambda_1 W_k}, \hat{\Lambda}),
			\end{equation}
			where $\lambda_1 \in \R$ and $W_k$ satisfies (\ref{W_k}) with $p>\frac{1}{2}$.
			\item\label{assump2(3)} There exists a constant $C>0$, and a non-decreasing function $h(\cdot)$ which is positive on any positive domain, such that for all $k < T_{\epsilon}$,
			\begin{equation}
			\E[ \Phi_{k+1} | \mathcal{F}_{k-1}^{M,F}] \leq \Phi_k - C h(\Lambda_k).
			\end{equation}
		\end{enumerate}
		
	\end{assumption}

	Based on Assumptions \ref{section2_assump2}, the following theorem (see \cite{STORM_nonconvex}) illustrates the upper bound on the expected number of iterations $T_{\epsilon}$ for obtaining an $\epsilon$-solution. 
	\begin{thm}\label{thm1}
		Let Assumption \ref{section2_assump2} holds. Then 	
		$$	\E [ T_{\epsilon} ]  \leq  \frac{p}{2p-1}( \frac{\Phi_{0}}{h(\hat{\Lambda})} + \frac{h(\Lambda_0)}{h(\hat{\Lambda})} + 1).$$			
	\end{thm}
	
	The analysis of the renewal-reward process in \cite{STORM_nonconvex} is appropriate for the stochastic process generated by Algorithm \ref{alg:1}. So here we omit the proof of Theorem \ref{thm1}. For more details, we refer the readers to Theorem 2.2 in \cite{STORM_nonconvex}. Theorem \ref{thm1} is very important to the following analysis of the complexity of Algorithm \ref{alg:1}. The difficulties lie in finding the non-decreasing function $h(\cdot)$ and the constant $\hat{\Lambda}$. The choice of trust-region radius indeed introduces some differences and difficulties, compared to the analysis of \cite{STORM_nonconvex} and \cite{S_line_search}. In the next part, we will give more analysis and discussions.

	%%%%%%%%%%%%%%%%%%%%%%%%%%%%%%%%%%%%%%%%%%%%%%%%%%%%%%%%%%%%%
	%%%%%%%%%%%%%%%%%%%%%%%%%%%%%%%%%%%%%%%%%%%%%%%%%%%%%%%%%%%%%	
	\section{Stochastic trust-region with probabilistic models and estimates}\label{sec:3}
	In this section, we propose a specific trust-region framework named STRME based on probabilistic models and estimates. The main steps are described as follows. 
	
	At each iteration $k$, given a current point $x_k$ and trust-region radius $\delta_k$, the model is built as 
	\begin{equation}
	m_k(x_k + d) = f_k + g_k^{T}d +\frac{1}{2}d^{T}B_kd, 
	\end{equation}
	to approximate $f(x)$ in $B(x_k, \delta_k)$. The quadratic model is simple and widely used in many trust-region algorithms. Of course, other models, for example the conic model (see \cite{Yuan_TR_review}), can also be applied to the framework as long as some requirements we stated are met. 
	
	The trust-region radius is defined as $\delta_k = \mu_k\left\|g_k \right\|$. Actually, one can try other more general choices that $\delta_k = \mu_{k}^{r_1}\left\|g_k \right\|^{r_2} $ with $r_1, r_2 \geq 0$ \cite{TR_Yuan_Trustregion}. For simplicity, we only consider the case that $r_1, r_2 = 1$. In the following steps, we choose to update the parameter $\mu_k$ (relative trust-region radius). 
	 Note that due to the randomness of the model, $\left\| g_k\right\| $ can be very small even zero even though the algorithm does not converge yet. In this case, it does not make any sense to continue the following process. Thus we add steps  in Algorithm \ref{alg:2} to check if $\left\|g_k \right\| > \epsilon $.
	
	The trial step $d_k$ is produced by minimizing the model $m_k(x_k + d)$ in a neighborhood of $x_k$ exactly or inexactly. Then we compute the random estimates $f_k^0$ and $f_k^d$ of $f(x_k)$ and $f(x_k + d_k)$ respectively to measure the actual function reduction. Once the trial step is obtained, we can use the ratio $\rho_k$, which is defined below, to judge how good the trial step $d_k$ is. Based on this criterion, if the trial point $x_{k} + d_k$ yields sufficient reduction, we accept the trial step $d_k$; otherwise, we reject it. At the end of each iteration, the relative trust-region $\mu_k$ is chosen according to the outcome of the iterates. The details of the algorithm are described as follows.

	\begin{algorithm}[H]
		\begin{algorithmic}[1]
			\STATE \textbf{Initialization}: Given an initial point $x_0$, $\gamma > 1$, $\eta_1\in(0,1)$, $\mu_0 \in (0, \mu_{\max})$ with $\mu_{\max} > 0$, $\epsilon=10^{-8}$; Set $k = 0$				
			%	\STATE  If $\left\|g_k \right\| \leq \epsilon $ then stop, else go to  step 3
			\STATE Construct a (random) model $m_k(x_k + d) = f_k + g_k^{T}d + \frac{1}{2}d^{T}B_kd $ that
			approximates $f(x)$ at $x_k$ with $d = x - x_k$
			\IF {$\left\|g_k \right\| \leq \epsilon $}
			\STATE return to step 2 until $g_k > \epsilon$
			\ENDIF
			\STATE Compute $\delta_{k} = \mu_{k}\left\|g_{k} \right\| $
			\STATE Compute the trial step $d_k = \arg\min_{\left\| d\right\| \leq \delta_k} m_k(x_k + d)$ such that $d_k$ satisfies Assumption \ref{sufficient_reduction}
			\STATE Obtain estimates $f_k^0$ and $f_k^d$ of $f(x_k)$ and $f(x_k + d_k)$
			\STATE Compute $\rho_k = \frac{f_k^0 - f_k^d}{m_k(x_k) - m_k(x_k + d_k)}$
			\IF{$ \rho_k \geq \eta_1$}
			\STATE$x_{k+1} = x_k + d_k$, $\mu_{k+1} = \min(\gamma\mu_k, \mu_{\max})$
			\ELSE 
			\STATE $x_{k+1} = x_k$, $\mu_{k+1} = \mu_k/\gamma$	
			\ENDIF	
			\STATE Set $k := k+1$, and go to step 2
		\end{algorithmic}
		\caption{Stochastic Trust-Region with Probabilistic Model and Estimates (STRME) }\label{alg:2}
	\end{algorithm}

	At each iteration, the trial step $d_k$ is computed to satisfy the well-known $ Cauchy\, decrease$ condition, which is given as follows. 	
	
	\begin{assumption}{\label{sufficient_reduction}}
		\begin{equation}
		m_k(x_k) - m_k(x_k + d_k) \geq \kappa_{fcd}\left\|g_k \right\|\min\left\lbrace \frac{\left\|g_k \right\| }{\left\| B_k\right\| }, \delta_k \right\rbrace. 
		\end{equation}
	\end{assumption}		

 Besides, in the case that $\delta_k = \mu_k\left\| g_k\right\|$, the condition $\left\| g_k\right\| \geq \eta_2 \delta_k $ in Algorithm \ref{alg:1} is equivalent to the condition $\mu_k \leq \frac{1}{\eta_2}$ . We notice that the parameter $\eta_2$ in algorithm \ref{alg:1} is usually very small. When $\eta_2$ is small, we can see that the condition $\mu_k \leq \frac{1}{\eta_2}$ is easily satisfied. We might as well let $\mu_{\max} \leq \frac{1}{\eta_2}$, thus the condition  $\left\| g_k\right\| \geq \eta_2 \delta_k $ in Algorithm \ref{alg:1} can be satisfied automatically.

	%%%%%%%%%%%%%%%%%%%%%%%%%%%%%%%%%%%%%%%%%%
	%%%%%%%%%%%%%%%%%%%%%%%%%%%%%%%%%%%%%%%%%%	
	\subsection{Theoretical properties of STRME }
	\label{sec:3:1}
	
	We are ready to present the theoretical properties of the framework described in Algorithm \ref{alg:2}. First, we give an assumption which states that the Hessian approximation matrix $B_k$ in model $m_k$ is uniformly upper bounded.

	\begin{assumption}{\label{hessian_bound}}
		There exists a constant $\kappa_{bhm} > 0$ such that, for all $k \geq 0$,
		\begin{equation*}
		\left\|B_k \right\| \leq \kappa_{bhm}.
		\end{equation*}	
	\end{assumption}
	
	We now provide some auxiliary lemmas to show that the decrease of the objective function $f(x)$ is guaranteed under some conditions. The following lemma states that if the model $m_k$ is true (fully linear) and the relative trust-region radius $\mu_k$ is upper bounded by a given number, then the actual reduction of the objective function is achieved. Although the theoretical content of Lemmas \ref{lem1} to \ref{lem_assump2(2)} closely relies on the existing arguments related to the STORM algorithm, for the  integrity of analysis, the proofs of these lemmas will still be attached in the Appendix.

	\begin{lem}{\label{lem1}}
		Suppose that model $m_k$ is true. If 	\vspace{-0.5em}
		\begin{equation*}
		\mu_k \leq \min\left\lbrace \frac{1}{\kappa_{bhm}}, \frac{\kappa_{fcd}}{8\kappa_{ef}}\right\rbrace,
		\end{equation*} 
		then \vspace{-0.5em}
		\begin{equation}
		f(x_k) -f(x_{k} + d_k)  \geq  \frac{\kappa_{fcd}}{4}\left\|g_k \right\|^2\mu_k. 
		\end{equation}			
	\end{lem}
	
	The next lemma states that the decrease of the objective function is achieved if the estimates $ \left\lbrace f_k^0, f_k^d \right\rbrace $ are  tight and iteration $k$ is successful.
	
	\begin{lem}{\label{lem2}}
		Suppose that estimates $f_k^0 $ and $ f_k^d $ are tight with $\epsilon_F \leq \frac{\eta_1\kappa_{fcd}}{8\mu_{\max}} $  and $\mu_k \leq \frac{1}{\kappa_{bhm}}$. If $d_k$ is accepted, then  
		\vspace{-0.5em}
		\begin{equation}
		f(x_k) - f(x_{k} + d_k) \geq \frac{\eta_1\kappa_{fcd}}{4}\left\|g_k \right\|^2\mu_k.
		\end{equation}	
	\end{lem}
	
	The following lemma shows when model $m_k$ and estimates $\left\lbrace f_k^0, f_k^d\right\rbrace $ are both sufficiently accurate, if $\mu_k$ is not too large, then the iteration will be successful.
	
	\begin{lem}\label{lem3}
		Suppose that model $ m_k $ is true,  and estimates $f_k^0 $ and $ f_k^d $ are tight with $\epsilon_F \leq \kappa_{ef} $. If \vspace{-0.5em}
		$$ \mu_k \leq \min \left\lbrace \frac{1}{\kappa_{bhm}}, \frac{\kappa_{fcd}(1-\eta_1)}{8\kappa_{ef}} \right\rbrace, $$
		\vspace{-0.5em}
		then  the k-th iteration is successful.	
	\end{lem}

	We now turn to consider the random process  $\left\lbrace \Phi_k, \Lambda_k\right\rbrace $ derived from the process generated from Algorithm \ref{alg:2}. The following analysis is based on the function
	\begin{equation}
	\Phi_k = \nu(f(X_k)- f^{\ast}) + (1-\nu)\frac{1}{L^2}\Lambda_k\left\|\nabla f(X_k) \right\|^2,
	\end{equation}
	for some  $\nu \in (0,1)$. It is obvious that $\Phi_k \geq 0$. Actually, the random variable $\Phi_k$ can be regarded as a kind of measure of progress to optimality. It plays an important role in the analysis of such trust-region algorithms. In Algorithm \ref{alg:2}, we update the trust-region radius as $\Delta_k = \Lambda_k \left\|G_k \right\| $. One may find that the trust-region radius $\Delta_k$ depends on the randomness introduced by the current model $M_k$, that is to say, $\Delta_k$ is not measurable with respect to $\mathcal{F}_{k-1}^{M,F}$. However, in STRME, we update the parameter $\Lambda_k$, which is completely determined by $\mathcal{F}_{k-1}^{M,F}$, instead of $\Delta_k$. Our approach avoids the difficulty of using the trust-region radius directly to analyze the complexity of STRME.  
	
	We have to admit that the theoretical analysis of STRME algorithm is similar to that of \cite{S_line_search} in spirit, but there are some distinctions between them. The most distinctive lies in the measure function $\Phi_k$. In \cite{S_line_search}, $\Phi_k$ consists of the three terms $f(X_k) - f^{\ast}$, $\frac{1}{L^2}\Lambda_k\left\|\nabla f(X_k) \right\|^2 $ and $\Delta_k^2$. However, in our analysis, only two of them, that is $f(X_k) - f^{\ast}$ and $\frac{1}{L^2}\Lambda_k\left\|\nabla f(X_k) \right\|^2 $, are used to evaluate the reduction of the algorithm. Actually, for STRME algorithm, trust-region $\Delta_k$ does not necessarily increase even if the iteration is successful.

	We aim to bound the expected number of iterations $\E[T_{\epsilon}]$. Before that, we have to prove Assumption \ref{section2_assump2} holds for the process $\left\lbrace \Phi_k, \Lambda_k\right\rbrace $. It is apparent that Assumption \ref{section2_assump2}(\ref{assump2(1)}) holds with the definition of $\Phi_k$ and $\mu_{k} \leq \mu_{\max}$ in Algorithm \ref{alg:2}. This assumption is not related with the convexity of the objective function, so it holds in all three cases we will consider later. Let us define the constant $\hat{\Lambda}$ in Assumption \ref{section2_assump2}(\ref{assump2(2)}) as follows:
	\begin{equation}\label{hat_Lambda}
	\hat{\Lambda} = \zeta, \text{where $\zeta$ is a constant such that} \,\, \zeta \leq \min\left\lbrace \mu_{\max},  \frac{\kappa_{fcd}(1-\eta_1)}{8\kappa_{ef}} \right\rbrace.
	\end{equation}
	In our analysis, we might as well claim that $\mu_{\max} \leq \min\left\lbrace \frac{1}{\kappa_{bhm}}, \frac{\kappa_{fcd}}{8\kappa_{ef}} \right\rbrace.$ For simplicity, we assume that $\Lambda_0 = \gamma^{i}\hat{\Lambda}$ and $\mu_{\max} = \gamma^{j}\hat{\Lambda}$ for some integers $i, j > 0$. As a result, for any $k > 0$, we have $\Lambda_k = \gamma^{i}\hat{\Lambda}$ for some integer $i$. Next, we will show that Assumption \ref{section2_assump2}(\ref{assump2(2)}) holds provided the constant $\hat{\Lambda}$ is defined as above. 
	
	\begin{lem}\label{lem_assump2(2)}
		Let $\alpha$ and $\beta$ safisfy that $\alpha\beta > \frac{1}{2}$, then we have Assumption \ref{section2_assump2}(\ref{assump2(2)}) holds with $W_k = 2(I_kJ_k - \frac{1}{2})$, $\lambda_1 = \log(\gamma)$, and p = $\alpha\beta$.
	\end{lem}
	
	Using the above Lemmas \ref{lem1} to \ref{lem3}, we can derive the following result.

	\begin{thm}\label{thm_Phi_decrease}
		Let Assumptions \ref{section2_assump1} hold with $\epsilon_F \leq \min \left\lbrace \frac{\eta_1\kappa_{fcd}}{8\mu_{\max}}, \kappa_{ef} \right\rbrace $ and $\kappa_f \leq 2\eta_1\kappa_{ef}$. Besides we assume that $f$ is $L$-smooth and Assumption \ref{hessian_bound} is satisfied. Then there exist a constant $\nu$ and sufficiently large $\alpha, \beta$
		satisfying the following conditions 
		\begin{equation}
		\frac{1-\nu}{\nu} \leq \min\left\lbrace  \frac{\kappa_{fcd}}{16\gamma \mu_{\max}^2}, \frac{\kappa_{fcd}L^2}{32\gamma(1+\kappa_{eg}\mu_{\max})^2}, \frac{\eta_1\kappa_{fcd}}{8\gamma\mu_{\max}^2}  \right\rbrace, 
		\end{equation}
		\begin{equation}
		\alpha\beta \geq \frac{4\gamma^2}{4\gamma^2+ (\gamma-1) },
		\end{equation}	
		such that 
		\begin{equation}
		\E[ \Phi_{k+1} - \Phi_k | \mathcal{F}_{k-1}^{M,F}] \leq - \frac{1}{2}\alpha\beta(1-\nu)(1-\frac{1}{\gamma})\frac{1}{L^2}\Lambda_k \left\|\nabla f(X_k) \right\|^2.
		\end{equation}
		
	\end{thm}

	\begin{rem}
		Due to the possibility of inaccurate model $M_k$ and estimates $\left\lbrace F_{k}^0, F_k^d\right\rbrace $, the function value of $f$ may increase.  So $\Phi_k$ is designed to balance the decrease and increase of $f(X_k)$. The above theorem shows that the decrease of expected $\Phi_k$ can be achieved by carefully choosing $\nu$ and the probability $\alpha$ and $\beta$.
	\end{rem}	
	
%%%%%%%%%%%%%%%%%%%%%%%%%%%%%%%%%%%%%%%%%%%%%%%%%%%%%%%%%%%%%%%%%	
	\subsection{Convergence rates for non-convex problems}
	\label{sec:3:2}
	We now show the global convergence rate of Algorithm \ref{alg:2} when $f$ is unknown to be convex, that is the following assumption holds.
	\begin{assumption}\label{L-Lipschitz}
		$f \in \mathcal{C}^{1}(\R^d) $ is bounded below by $f^{\ast}$ and {$L$-smooth}.	
	\end{assumption}
	Our goal is to bound the expected number of iterations until an $\epsilon$-solution occurs, i.e. $\E[T_{\epsilon}]$. The definition of $T_{\epsilon}$ is described as follows 
	\begin{equation}
	T_{\epsilon} = \inf\left\lbrace k\geq 0 : \left\|\nabla f(X_k) \right\| \leq \epsilon  \right\rbrace. 
	\end{equation}
	For all $k \leq T_{\epsilon}$, we know $\left\|\nabla f(X_k) \right\|\geq \epsilon $.
	Let us recall the definition of $\Phi_k$, i.e. $$\Phi_{k} = \nu (f(X_{k}) - f^{\ast}) + (1-\nu)\frac{1}{L^2}\Lambda_k\left\|\nabla f(X_k) \right\|^2.$$ In this case, $h(\cdot)$ can be defined as
	\begin{equation}
	h(\Lambda_k) =  C\Lambda_k \epsilon^2
	\end{equation}
	where $C  = \frac{1}{2}\alpha\beta(1-\nu)(1-\frac{1}{\gamma})$, which is non-decreasing on any positive domain. From Theorem \ref{thm_Phi_decrease}, we know that Assumption \ref{section2_assump2}(\ref{assump2(3)}) will hold if the conditions in Theorem \ref{thm_Phi_decrease} are satisfied.
	Applying the result in Lemma \ref{lem_assump2(2)}, Assumption \ref{section2_assump2}(\ref{assump2(2)}) holds under the conditions that $\alpha\beta > \frac{1}{2}$ and $\hat{\Lambda}$ is defined as (\ref{hat_Lambda}). Clearly Assumption \ref{section2_assump2}(\ref{assump2(1)}) holds.  Thus, we can conclude that Assumption \ref{section2_assump2} holds under certain conditions. Thus the conclusion of Theorem \ref{thm1} is true in this case. Then the following complexity result for Algorithm \ref{alg:2} can be achieved by simple analysis.
	\begin{thm}\label{thm_nonconvex}
		Suppose that $f$ satisfies Assumption \ref{L-Lipschitz}. Under the conditions in Theorem \ref{thm_Phi_decrease}, if $\alpha\beta > \frac{1}{2}$, then for Algorithm \ref{alg:2}, to achieve $\epsilon$-solution, the expected number of iterations is bounded as follows
		\begin{equation}\label{thm2_inequ}
		\E[T_{\epsilon}] \leq  \frac{\alpha\beta}{(2\alpha\beta -1)}(\frac{M}{\epsilon^2} + \mathcal{O}(1)\footnote{For simplicity, we use $\mathcal{O}(1)$ to denote the constant term of (\ref{thm2_inequ})}), 
		\end{equation}
		where $M = \frac{2\nu (f(x_0) - f^{\ast}) + 2(1-\nu)\frac{1}{L^2}\mu_0\left\|\nabla f(x_0) \right\|^2 }{\alpha\beta(1-\nu)(1-\frac{1}{\gamma})\hat{\Lambda}}$.
	\end{thm}
	\begin{rem}
		Note that the dependency of the expected number of iterations for obtaining an $\epsilon$-solution is $ \mathcal{O}(1/\epsilon^2)$, which is similar to that in \cite{STORM_nonconvex} for nonconvex problems.
	\end{rem}

%%%%%%%%%%%%%%%%%%%%%%%%%%%%%%%%%%%%%%%%%%%%%%%%%%%%%%%%%%%%%%	
	\subsection{Convergence rates for convex problems}
	\label{sec:3:3}
	In this part, we will analyze the expected complexity for STRME when $f$ is convex. First we give the following assumption.
	\begin{assumption}{\label{f_convex}}
		$f \in \mathcal{C}^{1}(\R^d) $ is convex. The level set $\mathcal{L} = \left\lbrace x\in \R^d : f(x) \leq f(x_0) \right\rbrace $ is bounded, and there exists a constant $ D > 0$ such that
		\begin{equation}
		\left\|x -x^{\ast} \right\| \leq D, \,\,\forall\, x \in \mathcal{L}.
		\end{equation} 
	\end{assumption}
	
	In convex setting, we aim to bound the expected $T_{\epsilon}$, which is defined as below
	\begin{equation}\label{stop_time}
	T_{\epsilon} = \inf\left\lbrace k\geq 0 : f(X_k) - f^{\ast} \leq \epsilon  \right\rbrace,
	\end{equation}
	for an $\epsilon$-solution. In this case, we define a function  
	\begin{equation}
	\Psi_{k} = \frac{1}{\nu\epsilon} - \frac{1}{\Phi_k},
	\end{equation}
	to replace $\Phi_k$ to measure the progress of the iterations. For $k \leq T_{\epsilon}$, $f(X_k) - f^{\ast} \geq \epsilon$, then we have $\Phi_k \geq \nu \epsilon $. Thus $\Psi_k \geq 0$, which is well-defined. In the later analysis, we can demonstrate that the random process $\left\lbrace \Lambda_k, \Psi_k\right\rbrace $ satisfies Assumption \ref{section2_assump2}. Therefore, Theorem \ref{thm1} can be applied to derive the upper bound of $\E[T_{\epsilon}]$. For more details, please refer to the proof of Theorem \ref{thm_convex} in the Appendix.

	\begin{thm}\label{thm_convex}
		We assume that $ f $ is $L$-smooth and satisfies Assumption \ref{f_convex}. Under the conditions in Theorem \ref{thm_Phi_decrease} and $\alpha\beta > \frac{1}{2}$, for Algorithm \ref{alg:2}, in order to achieve an $\epsilon$-solution, we have	
		\begin{equation}
		\E[T_{\epsilon}]  \leq  \frac{\alpha\beta}{(2\alpha\beta -1)}(\frac{M}{\epsilon} +  \mathcal{O}(1)),
		\end{equation}
		where $ M = \frac{2(\nu L  + (1-\nu)\mu_{\max} )^2D^2}{\alpha\beta\nu(1-\nu)(1-\frac{1}{\gamma})\hat{\Lambda}}.$ 	
		
	\end{thm}
	
	\begin{rem}
		
		The above theorem states that if $f$ is convex, the dependency on $\epsilon$ for the complexity bound is $\mathcal{O}(1/\epsilon)$, which is a stronger result than that in Theorem \ref{thm_nonconvex}. Note that the above result is the same as that of traditional trust-region for general convex problems. 
		
	\end{rem}
	
%%%%%%%%%%%%%%%%%%%%%%%%%%%%%%%%%%%%%%%%%%%%%%%%%%%%%%%%%%%%%%%%%%%%	
	\subsection{Convergence rates for strongly convex problems}
	\label{sec:3:4}
	In this part, we will derive the complexity bound for Algorithm \ref{alg:2} in the strongly convex setting. In this case, we assume the following assumption holds.

	\begin{assumption}{\label{f_stronglyconvex}}
		$f \in \mathcal{C}^{1}(\R^d) $ is strongly convex, i.e. there exists a constant $\sigma > 0$ such that 
		\begin{equation}
		f(x) \geq f(y) + \nabla f(y)^{T}(x-y) + \frac{\sigma}{2}\left\|x-y \right\|^2, \forall x, y \in \R^d.
		\end{equation}	
	\end{assumption}
	
	First we use the definition of $T_{\epsilon}$ as given in  (\ref{stop_time}):
	\begin{equation}
	T_{\epsilon} = \inf\left\lbrace k\geq 0 : f(X_k) - f^{\ast} \leq \epsilon  \right\rbrace.
	\end{equation}
	Our aim is to bound the expected number of iteration i.e. $\E[T_{\epsilon}]$, to obtain an $\epsilon$-solution for strongly convex problems. Instead of using $\Phi_k$ to measure the progress towards optimality, we define a function 
	\begin{equation}
	\Psi_{k} = \log(\Phi_k) + \log(\frac{1}{\nu\epsilon}).
	\end{equation}
	For $k \leq T_{\epsilon}$, $f(X_k) - f^{\ast} \geq \epsilon$, then we have $\Phi_k \geq \nu \epsilon $, which implies that $\Psi_k \geq 0$. So the definition of $\Psi_k$ is reasonable. In the later part, we will obtain the upper bound of $\E[T_{\epsilon}]$ with the help of random process $\left\lbrace \Lambda_k, \Psi_k\right\rbrace $. The details will be shown in the proof of the following theorem.

	\begin{thm}\label{thm_stronglyconvex}
		Assume that $ f $ is $L$-smooth and satisfies Assumption \ref{f_stronglyconvex}. If the conditions in Theorem \ref{thm_Phi_decrease} and $\alpha\beta > \frac{1}{2}$ hold, then for Algorithm \ref{alg:2}, in order to achieve an $\epsilon$-solution, we have	
		
		\begin{equation}
		\E[T_{\epsilon}]  \leq \frac{\alpha\beta}{2\alpha\beta -1} (M\log(\frac{1}{\epsilon}) + \mathcal{O}(1)),
		\end{equation}
		where $ M = \frac{(\frac{\nu}{2\sigma}  + (1-\nu)\frac{1}{L^2}\mu_{\max} )}{\alpha\beta(1-\nu)(1-\frac{1}{\gamma})\hat{\Lambda}}.$
	\end{thm}
	Similarly, the proofs of the above theorem will be given in the Appendix.
	\begin{rem}
		The result in Theorem \ref{thm_stronglyconvex} show Algorithm \ref{alg:2} takes at most $\mathcal{O}(\log(1/\epsilon))$ iterations in expectation, to achieve an $\epsilon$-solution. The result  coincides with that for trust-region methods in deterministic setting.
		
	\end{rem}
	
	%%%%%%%%%%%%%%%%%%%%%%%%%%%%%%%%%%%%%%%%%%%%%%%%%%%%%%%%%%%%		
	%%%%%%%%%%%%%%%%%%%%%%%%%%%%%%%%%%%%%%%%%%%%%%%%%%%%%%%%%%%%
	%%%%%%%%%%%%%%%%%% Numerical Experiments %%%%%%%%%%%%%%%%%%%%
	%%%%%%%%%%%%%%%%%%%%%%%%%%%%%%%%%%%%%%%%%%%%%%%%%%%%%%%
	\section{Numerical experiments}\label{sec:4}
	In this section, we empirically test our STRME algorithm and compare its performance with STORM\cite{STORM} and some related algorithms. 
	
	We test on two type of problems: (i) regularized logistic regression problem, which is strongly convex; (ii) deep neural networks, which is highly non-linear and non-convex. The function value of training data (called training loss) and accuracy (percentage of correctly classified testing data) are adopted as criteria to measure the performance of all the algorithms that are tested. For all those algorithms, we compare these criteria against the number of effective pass through the data, that is total gradient calls divided by $N$(training data size). All algorithms were terminated when the maximum budget of the gradient evaluations is larger than the  maximum value $SFO_{\max}$ we have set.

	All algorithms are implemented in Anaconda3 (python 3.6.2) under Windows 7 operating system on Dell desktop with Intel(R) Core(TM) i7-4790U CPU @3.6GHz, 8GB Memory.
	\subsection{How to implement the probabilistic models and estimates to satisfy Assumption \ref{section2_assump1}}\label{sec:4:1}
 In these numerical experiments, we focus on the derivative-based problems where $f(x,\omega)$ and $\nabla f(x, \omega)$ are available. At each iteration point $x_k$, we assume that the noise $\omega$ in function value and gradient computation is conditional unbiased and the corresponding variance is conditional bounded for all $f$, i.e.	
	\begin{equation}
	\begin{split}
	\E [f(x_k;\omega) \,|\, \mathcal{F}_{k-1}^{M,F} ] = f(x_k)  \quad & \text{and}\quad  Var[f(x_k;\omega)  \,|\, \mathcal{F}_{k-1}^{M,F}] \leq  V_f ;\\	 
	\E [\nabla f(x_k;\omega)  \,|\, \mathcal{F}_{k-1}^{M,F}]= \nabla f(x_k) \quad  & \text{and} \quad   Var[\nabla f(x_k;\omega)  \,|\, \mathcal{F}_{k-1}^{M,F}] \leq  V_g.	
	\end{split}
	\end{equation}

	We now discuss how to obtain $\alpha$-probabilistically $\kappa$-fully linear models and $\beta$-probabilistically $\epsilon_{F}$-accurate estimates. In the later analysis, we employ the standard sample averaging approximation technique \cite{Sample_average} to construct the sufficient accurate model and estimates. Let
	\begin{equation}\label{prob_model}
	f_k = \frac{1}{|S_k|}\sum_{i\in S_k} f(x_k;\omega_i), \quad \text{and}\,\,g_k = \frac{1}{|S_k|}\sum_{i\in S_k} \nabla f(x_k;\omega_i),
	\end{equation}
	where $\omega_i$ are the i.i.d.
	and finite realizations of the noise $\omega$ and the sample set $S_k\subseteq \left\lbrace 1,2,\cdots, N \right\rbrace $ of size $\left|S_k \right|=p_k $.
	The local approximation model can be constructed as $m_k(x) = f_k + g_k^{T}(x-x_k) + \frac{1}{2}(x-x_k)^{T}B_k(x-x_k)$, where $x \in B(x_k, \delta_k)$.
	We recall the definitions in Section \ref{sec:2:1}. In order to satisfy Assumption \ref{section2_assump1}(\ref{assump1(1)}),we need the following conditions hold
	\begin{equation}\label{ineq_fg}
	\begin{split}
	P(|f(x) - m_k(x)| \geq \kappa_{ef}\delta_k^2\,\, | \,\, \mathcal{F}_{k-1}^{M,F}) & \leq 1- \alpha'; \\ 
	P(\left\|\nabla f(x) -\nabla m_k(x) \right\|  \geq \kappa_{ef}\delta_k \,\, | \,\, \mathcal{F}_{k-1}^{M,F})  &\leq 1- \alpha' (\alpha = \alpha'^2)
	\end{split}
	\end{equation}
	for all $x \in B(x_k, \delta_k)$.
	By Chebyshev's inequality in Lemma \ref{Chebyshev_inequ}, at current point $x_k$, for any $v > 0$, we have 
	\begin{equation}\label{ineq_cheb}
	P[\left|f_k - f(x_k) \right| \geq  v \,\, | \,\, \mathcal{F}_{k-1}^{M,F}] \leq \frac{V_f}{p_kv^2},\quad  \text{and} \,\, P[\left\|g_k - \nabla f(x_k) \right\|\geq  v \,\, | \,\, \mathcal{F}_{k-1}^{M,F} ] \leq \frac{V_g}{p_kv^2}.
	\end{equation} 
By (\ref{ineq_cheb}), the conditions (\ref{ineq_fg}) hold at $x=x_k$,
as long as the sample rate $p_k$ satisfies the following condition that
	\begin{equation} \label{sample_g}
	p_k \geq \max\left\lbrace \frac{V_f}{(1-\alpha')\kappa_{ef}^2\delta_k^4}, \frac{V_g}{(1-\alpha')\kappa_{eg}^2\delta_k^2} \right\rbrace.
	\end{equation}
	For $ \forall x\in B(x_k, \delta_k)/\left\lbrace x_k\right\rbrace $, if $\left\|g_k - \nabla f(x_k) \right\| \leq \kappa_{eg} \delta_k$, we have 
	\begin{equation}
	\begin{split}
	& \left\|\nabla m_k(x) - \nabla f(x) \right\|  \\
	= & \left\|\nabla g_k + B_k(x-x_k) - \nabla f(x) \right\|\\
	 =  & \left\|\nabla g_k + B_k(x-x_k) - \nabla f(x_k) + \nabla f(x_k)-\nabla f(x) \right\| \\
	 \leq & \left\|g_k - \nabla f(x_k) \right\|+ L\left\|x_k -x \right\|+\left\|B \right\|\left\|x-x_k \right\|     \\	
	 \leq & (\kappa_{eg} + L + \kappa_{bhm} )\delta_k.
		\end{split}
	\end{equation}
	The second inequality relies on the fact that $\left\|x-x_k \right\|\leq \delta_k $, $f$ is $L-$smooth, and the second-order matrix $B_k$ satisfies Assumption \ref{hessian_bound}. 
	Thus we can conclude that if the sample rate $p_k$ satisfies (\ref{sample_g}), condition (\ref{ineq_fg}) will hold for any $x\in B(x_k, \delta_k)$. Consequently, $m_k(x)$ is $\alpha$-probabilistically $\kappa$-fully linear models. In practice, the function value is not explicitly computed at model $m_k$, so we only require that $ p_k \geq \frac{V_g}{(1-\alpha')\kappa_{eg}^2\delta_k^2}.$
	
	Similarly, to obtain $\beta$-probabilistically $\epsilon_{F}$-accurate estimates, let estimates $\left\lbrace f_k^0, f_k^d \right\rbrace $ be
	\begin{equation}
	f_k^0 = \frac{1}{q_k}\sum_{i\in S_k^0} f(x_k; \omega_i),\quad \text{and}\,\,  f_k^d = \frac{1}{q_k}\sum_{i\in S_k^0} f(x_k + d_k; \omega_i),	
	\end{equation} 
	where $S_k^0 \subseteq \left\lbrace 1,2,\cdots, N \right\rbrace $ and $\left|S_k^0 \right|=q_k $.	In order to satisfy Assumption \ref{section2_assump1}(\ref{assump1(2)}) such that	
	\begin{equation}
	\begin{split}
	P(|f(x_k) - f_k^0| \geq \epsilon_{F}
	\delta_k^2  \,\, | \,\, \mathcal{F}_{k-1/2}^{M,F}) \leq 1- \beta, &\, \text{and}\,\,
	P( |f(x_k + d_k) - f_k^{d} | \geq \epsilon_F \delta_k^2 \,\, | \,\, \mathcal{F}_{k-1/2}^{M,F}) \leq 1- \beta,
	\end{split}
	\end{equation}		
	we require 
	\begin{equation}\label{sample_f}
	q_k \geq \frac{V_f}{(1-\beta)\epsilon_{F}^2\delta_k^4} (\approx  \mathcal{O}(\frac{1}{\delta_k^4})).
	\end{equation}	
	By the H\"{o}lder's inequality for expectation that $\E [ ab ] \leq (\E[a^2])^{1/2}(\E[b^2])^{1/2} $ for $a,b > 0$, we have 
	\begin{equation}
	\E[\frac{\left| f_k^0 - f(x_k)\right| }{\epsilon_{F}\delta_k^2} | \mathcal{F}_{k-1/2}^{M,F} ] \leq (\E[1 |\mathcal{F}_k^{M,F} ])^{\frac{1}{2}}(\E[\frac{\left| f_k^0 - f(x_k)\right|^2 }{\epsilon_{F}^2\delta_k^4} | \mathcal{F}_{k-1/2}^{M,F}])^{\frac{1}{2}}. 
	\end{equation}
	Since $q_k \geq \frac{V_f}{(1-\beta)\epsilon_{F}^2\delta_k^4}$, we have $\E[\frac{\left| f_k^0 - f(x_k)\right|^2 }{\epsilon_{F}^2\delta_k^4} | \mathcal{F}_{k-1/2}^{M,F}]  \leq  \frac{V_f}{q_k\epsilon_{F}^2\delta_k^4}\leq \mathcal{O}(1) $. Thus there exists a constant $\kappa_{f} > 0$ such that $\E[\left| f_k^0 - f(x_k)\right| | \mathcal{F}_{k-1/2}^{M,F}]$ is bounded by $\kappa_{f}\delta_k^2$. In the same way, similar result can be applied to $f_k^d - f(x_k + d_k)$. We can claim that Assumption \ref{section2_assump1}(\ref{assump1(3)}) holds.
	Thus, we have shown how to construct the random model $m_k$ and estimates $\left\lbrace f_k^0, f_k^d \right\rbrace $ to satisfy Assumption \ref{section2_assump1}.	
	
	However, we can see that the trust-region radius $\delta_k$ is unknown before the sample rate is chosen as in (\ref{sample_g}) and (\ref{sample_f}). It brings a challenge to construct a fully linear model in practice. First of all, obviously, in the case that the total sample set is large but limited, such sample rate $p_k$ must exist, such as $p_k = N$, to make the model fully linear. But it is really difficult to give a computable condition of $p_k$ in theory, and prove that under such condition the model is fully linear with probability. Thus we settle for the minimal sample rate of increase in $p_k$ as \cite{byrd2012sample}.
	As previous discussed in (\ref{sample_g}), the sample size $p_k$ satisfies 
	\begin{equation}\label{p_condition}
	p_k\geq \frac{V_g}{(1-\alpha')\kappa_{eg}^2\delta_k^2} \geq  \frac{V_g}{(1-\alpha')\kappa_{eg}^2\mu_{\max}^2\left\| g_k\right\|^2} \approx \mathcal{O}(\frac{V_g}{\left\| g_k\right\|^2}).
	\end{equation}
	In order to guarantee the fully linear property, from the analysis in Section \ref{sec:3}, the critical is that there exists a constant $G_1 > 0$ such that the following condition holds
	\begin{equation}\label{g_F}
	\left\|g_k - \nabla F(x_k) \right\| \leq G_1\left\| g_k\right\|.
	\end{equation} 
   The above condition yields 
   \begin{equation}\label{g_F2}
   \left\|g_k \right\| \geq  \frac{1}{G_1 +  1}\left\|\nabla F(x_k) \right\|.
   \end{equation}
	In \cite{byrd2012sample}, under the condition (\ref{g_F}), they have analyzed gradient-based mini-batch optimization algorithm for the strongly convex problem that the sample size should grow geometrically with iteration $k$. Intuitively, it make sense because if condition (\ref{g_F}) holds, $\left\|\nabla F(x_k)\right\|^2 $ is geometrically decreasing in this case, so does $\left\|g_k \right\|^2 $. Theorem \ref{thm_stronglyconvex} shows that our algorithm achieves the same complexity rate as Theorem 4.1 in \cite{byrd2012sample}. By (\ref{p_condition}), in the strongly convex case, we impose the sample rate $p_k$ to be exponentially increased, i.e. $p = a^{k}$ for some $a > 1$. However, for the non-convex problem, the situation is different.
	From the analysis in Section \ref{sec:3:2}, to obtain an $\epsilon$-solution ($\left\|\nabla F(x_k) \right\| \leq \epsilon $), we have $T_{\epsilon} \approx \mathcal{O}(\frac{1}{\epsilon^2})$. From (\ref{p_condition}) and (\ref{g_F2}), the sample rate $p_k$ should be linearly increasing with iteration $k$.
	
	As we know, the situation in practice is very complicated. Thus we require that the sample rate grows exponentially with iteration $k$ for regularized logistic regression problem and linearly increases with iteration $k$ for the deep neural networks, respectively.
 
	%%% 目前是这样，就是对于强凸的问题，p可以按照指数增加，所以这里的数值实验还是要加以修正，但是因为STORM数值实验中用的是线性增加，所以我们需要做线性和指数增加的对比。但是对于非凸的问题，我们可以分析p只需要是线性增加就可以了。

	\subsection{Experimental results on regularized logistic regression problem}\label{sec:4:2}
	In this subsection, 
	we consider the following smooth (strongly convex) regularized logistic loss problem considered in \cite{STORM}:
	
	\begin{equation}\label{Fun_numerical}
	F(x) = \frac{1}{N}\sum_{i=1}^{N}\log(1 + \exp(-b_i(a_i^{T}x ))) + \frac{\lambda}{2} \left\|x \right\|^2, x \in \R^d,
	\end{equation}
	where $\left\lbrace (a_i, b_i)\right\rbrace_{i=1}^{N} $ is a training sample set with $a_i\in\mathbb{R}^d$ being the feature vector and  $b_i \in \left\lbrace -1, +1 \right\rbrace $ being the corresponding label. And $\lambda \geq 0$ is the regularization parameter. As in the typical machine learning setting, we assume that $n$ is very large and $N \gg d$. So computing $F(x)$, as well as $\nabla F(x)$ and $\nabla^2 F(x)$ are very expensive. In our work, we randomly (without replacement) choose a subset $I_k \subseteq \left\lbrace 1,2,\cdots, N \right\rbrace $ to estimate the quantities in our algorithms. For the algorithms that need to compute Hessian matrix, the same sample is drawn for gradient and Hessian evaluations. In this setting, we re-sample the sample set for $f_k^0$ and $f_k^d$.
	
	We compare our algorithm STRME with STORM which is implemented as algorithm 5 in \cite{STORM}. In our numerical experiments, we construct two versions of STRME: one which only computes the stochastic gradients and sets $B_k = 0$ is the first-order version (called STRME-1st), the other one which in addition to the stochastic gradients, computes stochastic Hessian estimators is the second-order version (called STRME-2st). We use the classic dogleg method in \cite{Numerical-Op} to solve the second order subproblem , and the corresponding algorithms, we call STRME-dogleg and STORM-dogleg. For the details of the implementation, one can refer to Algorithm \ref{alg:3} in Appendix. Besides, We compare against a special adaptive solver AdaGrad \cite{AdaGrad}, which takes the adaptive step size but does not have to compute function value, and only computes the average stochastic gradients.
	
	For problem (\ref{Fun_numerical}), all algorithms were tested with different input parameters. We set $x_0 = 0$ as the starting point for all algorithms. For the three quantities gradient, Hessian and ($f_k^0$, $f_k^d$) evaluations, we adopt the linearly increased sample rule that $b_k = \min \left\lbrace b_{\max}, \max\left\lbrace t_0 k + b_0, \frac{1}{\delta^2} \right\rbrace \right\rbrace $, where $b_0=d+1$, $b_{\max} = N$, $t_0=100$ for STORM and our STRME algorithm as in \cite{STORM}.  For AdaGrad, the mini-batch size $b = d+1.$ We set the same random seed to generate random sample sequence for all the algorithms.
		
	The regularized logistic regression problem we consider in this subsection is strongly convex. As we discussed in section \ref{sec:4:1}, to make sure that the model is fully linear, the sample size should be exponentially increased. Thus we also tested the mini-batch size $b_k = \min \left\lbrace b_{\max}, \max\left\lbrace b_0 a^{k}, \frac{1}{\delta^2} \right\rbrace \right\rbrace $, where $b_0=d+1$, $b_{\max} = N$, $t_0=100$, and $ a\in \left\lbrace 1.1, 1.3, 1.5, 2\right\rbrace $ for our algorithm STRME and STORM. And we compare the two sample rules as well at the end of this part.

	In this subsection, we test on two datasets a9a and ijcnn1 from the LIBSVM website \footnote{\url{https://www.csie.ntu.edu.tw/ cjlin/libsvmtools/datasets/}}.
	We list the datasets in Table \ref{tab:1}, in which $n$ denotes the total sample size,  and $d$ is the dimension of the dataset, and $\lambda$ is the regularization parameter.  We use 0.95 partition of the data as the training set, and the remaining as the testing set, just like in \cite{STORM}.

	\begin{table}[H]
	\centering
		\caption{Datasets for regularized logisitic regression}
		\label{tab:1}      % Give a unique label
		\begin{tabular}{llll}
			\hline\noalign{\smallskip}
			dataset & $n$ & $d$ & $\lambda$\\
			\noalign{\smallskip}\hline\noalign{\smallskip}
			a9a & 32561 & 123 & $10^{-4}$ \\
			ijcnn1 & 49990  & 22 & $10^{-4}$\\
			\noalign{\smallskip}\hline
		\end{tabular}
	\end{table}

	\begin{figure}[H]
		\centering	
		\includegraphics[width=0.75\textwidth]{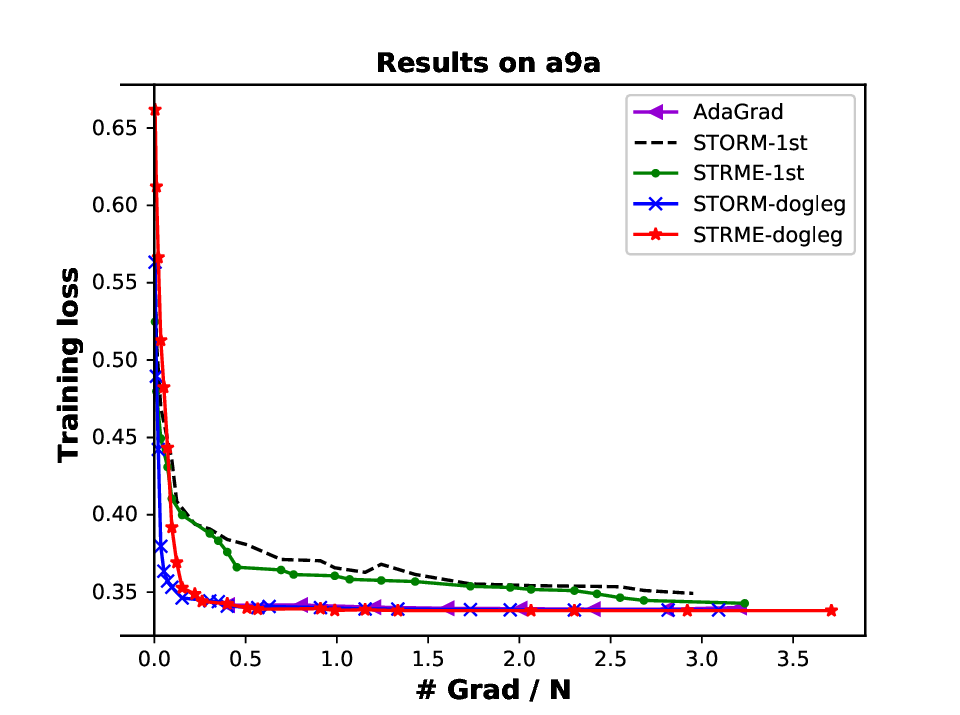}
		\hfill 
		\includegraphics[width=0.75\textwidth]{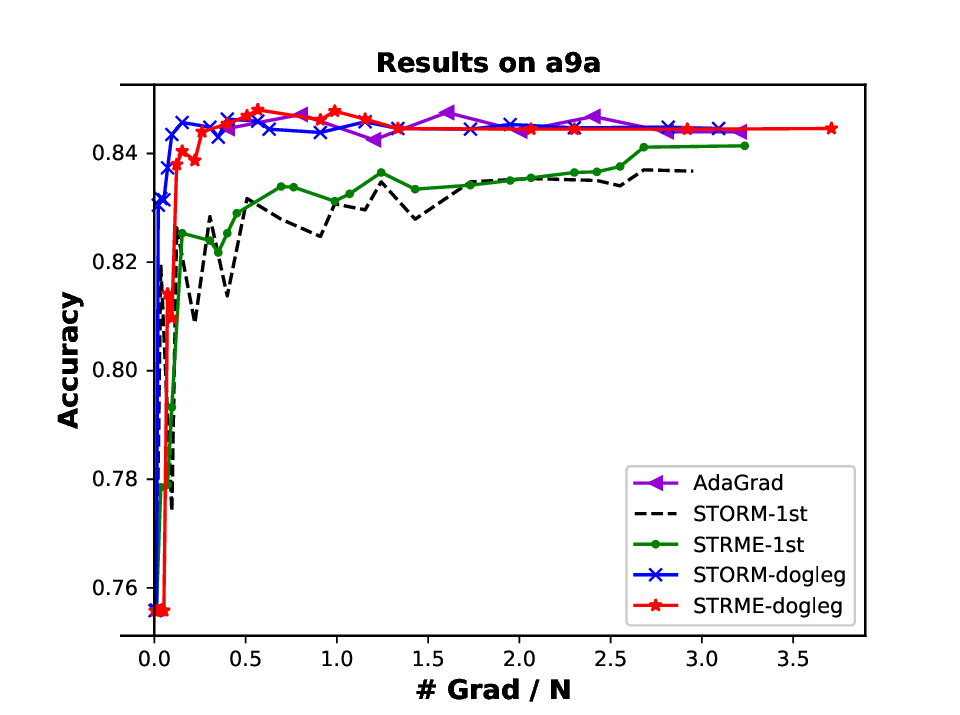}	
		\caption{Training regularized logistic regression on a9a}
		\label{fig:1}       % Give a unique label
	\end{figure}

	\begin{figure}[H]
		\centering	
		\includegraphics[width=0.75\textwidth]{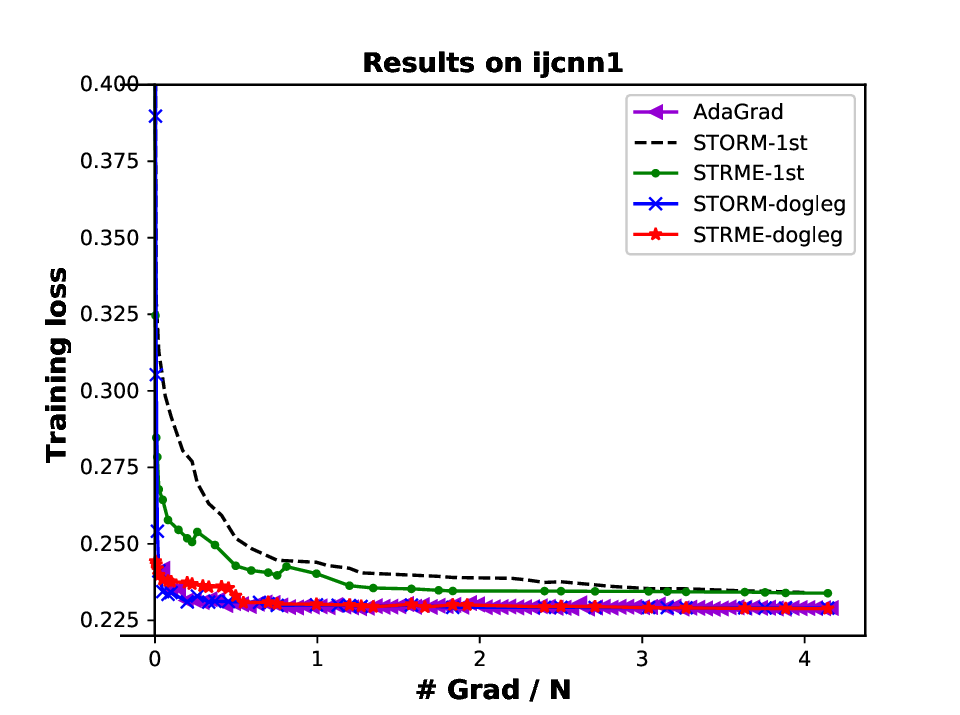}
		\includegraphics[width=0.75\textwidth]{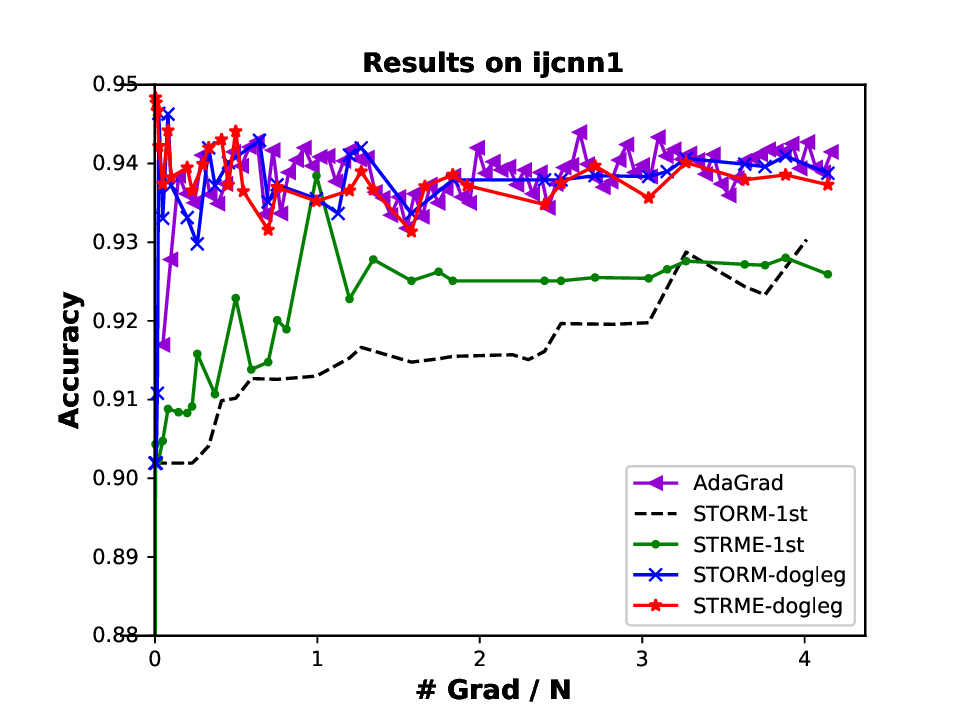}
		\caption{Training regularized logistic regression on ijcnn1}
		\label{fig:2}
	\end{figure}
	In Figure \ref{fig:1}, we present the results on the a9a dataset.
	For STORM-1st and STORM-2st, we use the same parameters in \cite{STORM}: $\delta_{\max} =10, \delta_0 =1, \gamma=2, \eta_1=0.1, \eta_2=0.001$. For STRME-1st and STRME-2st, the following parameters are used: $\mu_{\max} = 10^{3}, \mu_0 =1,  \eta_1 = 0.1, \gamma=2$. We set step size $\eta = 1$ for AdaGrad. 
	
	In Figure \ref{fig:2}, we report the results on the ijcnn1 dataset. The parameters $\mu_0 =10, \mu_{\max} = 10^{3}, \eta_1 = 0.1, \gamma=2$ are set for STRME-1st and $\mu_0 = 1$, others are the same for STRME-dogleg. For STORM-1st and STORM-dogleg, we set $\delta_0 = 1, \delta_{\max} = 10, \eta_1 = 0.1, \gamma = 2$.  For AdaGrad, we set step size $\eta = 1$.
	
	From Figure \ref{fig:1} and \ref{fig:2}, we can find that the proposed STRME is comparable to STORM in this setting, both in terms of the training function value and accuracy.
		
	At the end of this part, we compare the exponentially increased sample rule with the linearly increased case in Figure \ref{fig:3}. We choose the exponential ratio $a=1.1$ from $\left\lbrace 1.1, 1.3, 1.5, 2\right\rbrace $ for all the experiments. The other parameters are the same as the linearly increased sample. For the sake of distinction, we use STRME-L and STRME-E to denote the linearly increased and exponentially increased sample rules respectively, so does STORM. As a whole, from Figure \ref{fig:3}, we can see that the performance on the two sample rules is very similar.

		\begin{figure}[H]
			\centering	
			\includegraphics[width=0.75\textwidth]{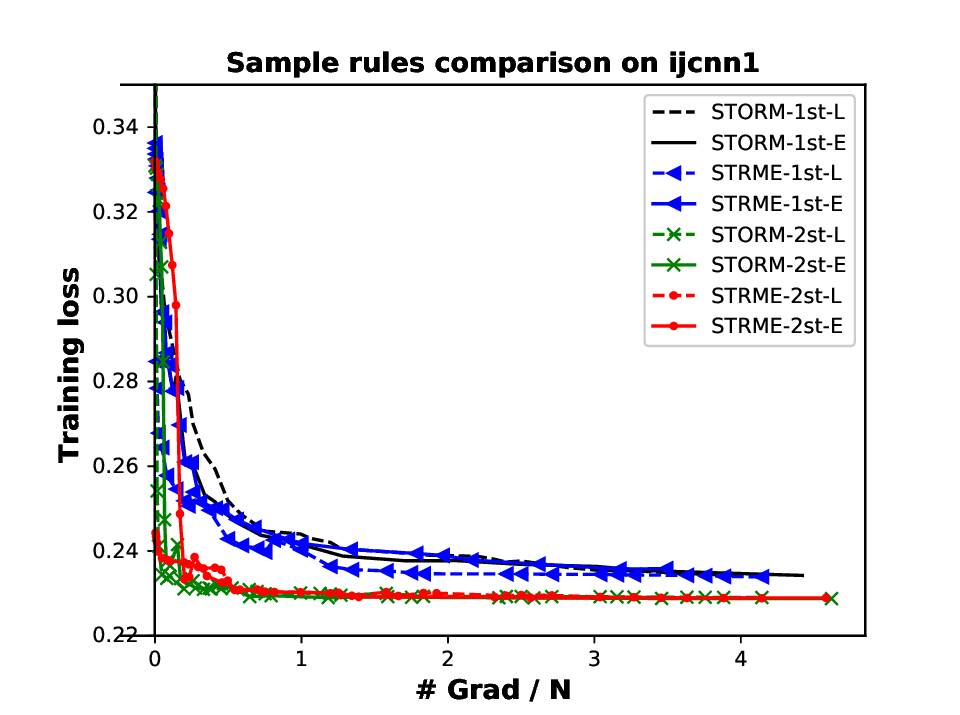}
			\includegraphics[width=0.75\textwidth]{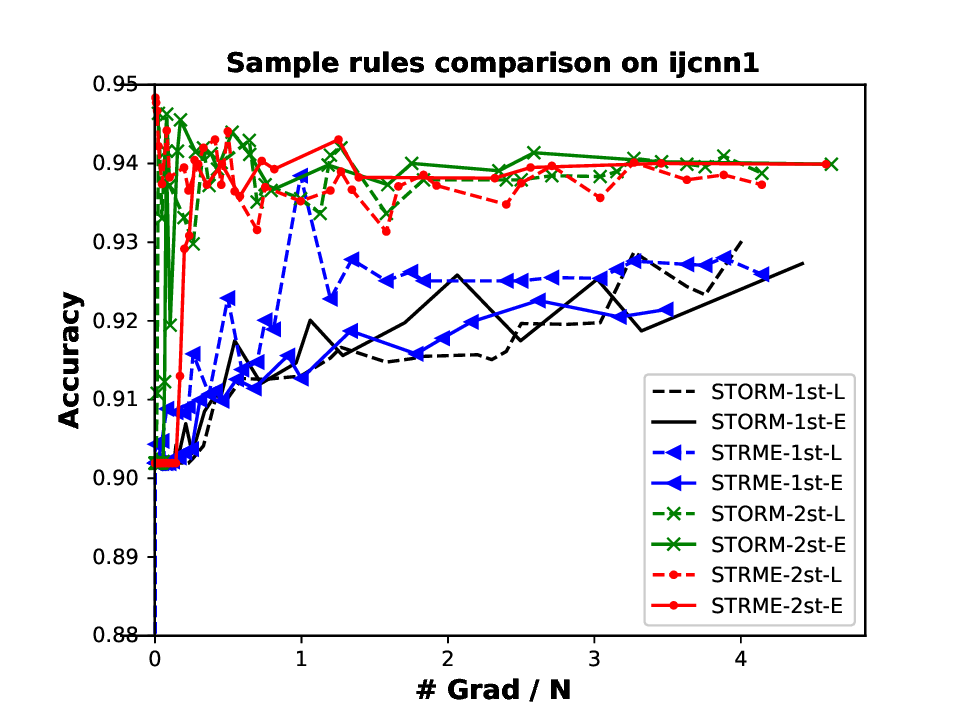}
			\caption{ Sample rules comparison on ijcnn1}
			\label{fig:3}
		\end{figure}

	\subsection{Experimental results on a simple deep neural network(DNN)}\label{sec:4:3}	
	
	In this subsection, we
	consider to train a fully-connected 2-layer net with 50 hidden units (784-50-10) neural networks with MNIST\footnote{\url{http://deeplearning.net/data/mnist/}}, a benchmark dataset of  handwritten digits. We used softmax output, sigmoid hidden functions, and the cross-entropy error function. The $l_2$ regularization parameter $\lambda = 10^{-3}$, suggested in \cite{SVRG-nonconvex}.  
	
	We compare our algorithm STRME with STORM\cite{STORM}. As in the previous subsection, we construct two versions of STRME: one which only computes the stochastic gradients and set $B_k = 0$ is the first-order version, the other one which in addition to the stochastic gradient, computes the quasi-Newton matrix $B_k$ to approximate the true Hessian matrix is the second-order version.
	Besides, we make an implementation of the SdLBFGS \cite{SdLBFGS}, which is an efficient second-order algorithm for the non-convex problem, to compare with STRME and STORM.	
	
	In our work, we run one epoch mini-batch SGD algorithm to obtain an initial point for all algorithms. In our implementation, we randomly (without replacement) choose the subset $I_k \subseteq \left\lbrace 1,2,\cdots, n \right\rbrace $ to estimate the gradient and Hessian pair $(g_k, B_k)$ and objective function value pair ($f_k^0$, $f_k^d$). We have attempted the three cases: (i) sample gradient and Hessian pair, and re-sample $f_k^0$ and $f_k^d$  independently; (ii) sample the two pair independently; (iii) sample the two pairs with the same subset. However, for the first two cases, the results are not satisfactory. Therefore, in our numerical experiments, we only consider the last one. For STORM and STRME, we set the mini-batch size $b_k = \min \left\lbrace b_{\max}, \max\left\lbrace t_0 k + b_0, \frac{1}{\delta^2} \right\rbrace \right\rbrace $, where $b_0=d+1$, $b_{\max} = N$, $t_0=10$. 
	
	\subsubsection{Experimental results on the first-order probabilistic model}\label{sec:4:3:1}
	
	In this part, we first construct a first-order probabilistic model $m_k$ at $x_k$, i.e. 
	\begin{equation}
	m_k(x_k + d) = f_k + g_k^{T}d, 
	\end{equation}
	where $f_k$ and $g_k$ are computed as (\ref{prob_model}), to test  the STRME (STRME-1st) framework. Actually, in practice, we do not need to compute $f_k$. In our numerical experiments, we compare the numerical performace of STRME-1st with the related STORM-1st.

	We now give details of parameters in the proposed STRME-1st and STORM-1st. For STRME-1st, $\mu_0$ and $\mu_{\max}$ are two important parameters.  The parameters $\mu_0$ is chosen from  $\left\lbrace 0.01, 0.1, 1\right\rbrace $, and the best $\mu_0$ is achieved at $\mu_0 = 0.1$. Compared to $\mu_0$, the parameter $\mu_{\max} $ is more important in the later iteration process. Therefore the range of $\mu_{\max}$ is more elaborate. In our implementation, let $\mu_{\max} \in \left\lbrace 1/2, 1, 2, 2^2, 2^3 \right\rbrace $, and the best tuned $\mu_{\max}$ is achieved at $\mu_{\max} = 2$. For STORM-1st, we test on $\delta_0 \in \left\lbrace 0.01, 0.1, 1\right\rbrace $ and $\delta_{\max} \in \left\lbrace 1/2, 1, 2, 2^2, 2^3 \right\rbrace $, the best performance is achieved with the inputs $(\delta_0, \delta_{\max}) = (0.1, 1)$.  The results are showing in the Figure \ref{fig:4}.

	\begin{figure}[H]
		\centering		
		\includegraphics[width=0.75\textwidth]{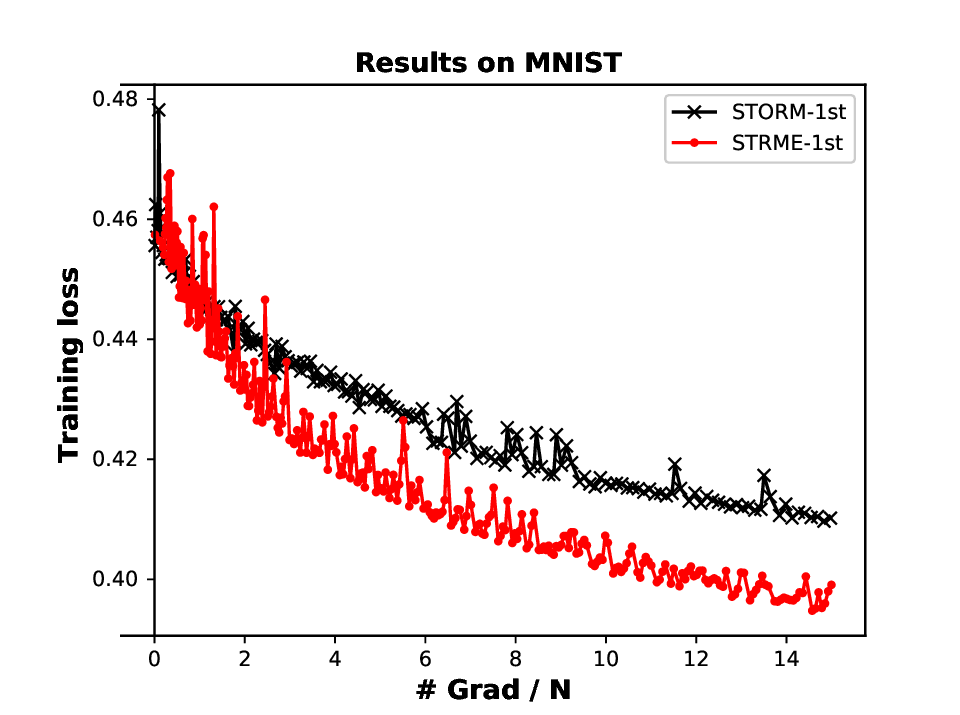}
		\includegraphics[width=0.75\textwidth]{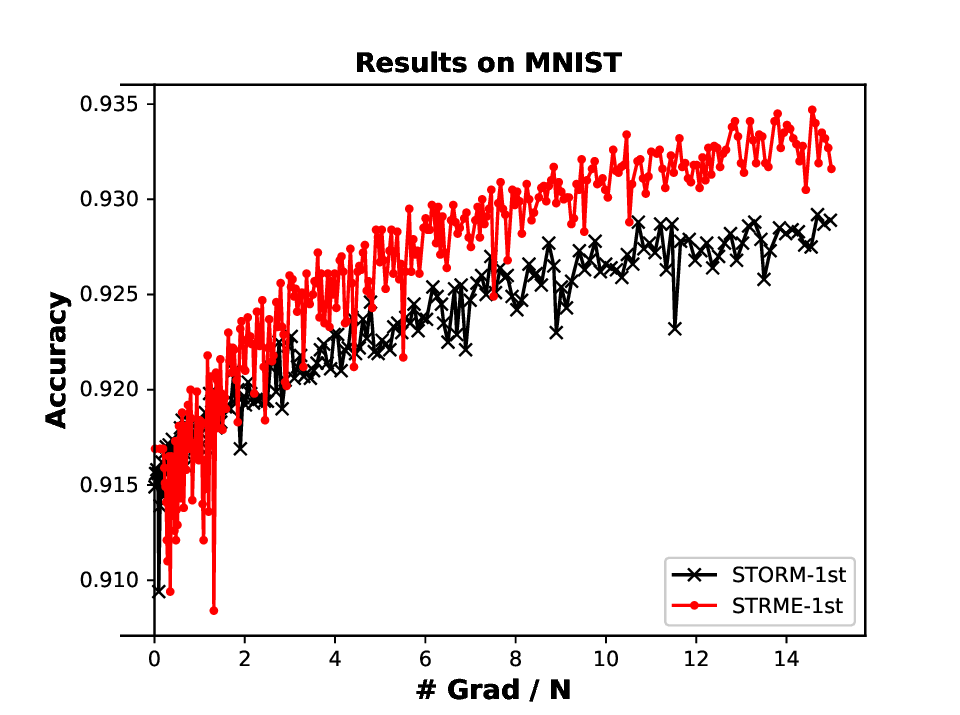}
		\caption{Training DNN on first-order probabilistic model with $\mu_0 =0.1, \mu_{\max} = 2, \eta_1 = 0.1, \gamma=2$ for STRME-1st. For STORM-1st, we set $\delta_0 = 0.1, \delta_{\max} = 1, \eta_1 = 0.1, \gamma = 2$.  }
		\label{fig:4}
	\end{figure}
	
		From Figure \ref{fig:4}, we can see that STRME-1st is better than STORM-1st with the model constructed as the beginning of this section. Not surprisingly, the trust-region radius constructed by STRME can make better use of gradient information. %Beyond that it can be seen that STRME-1st is comparable with AdaGrad. 
		However, we can construct a more efficient second-order model to make the STRME algorithm perform better.

	\subsubsection{Experiments on second-order probabilistic models}\label{sec:4:3:2}
	In this part, we construct a specific second-order probabilistic model and implement STRME algorithm framework in deep neural network. In our work, we consider the limited memory symmetric  rank one method (L-SR1) to generate the second order quasi-Newton matrix $B_k$, and build up the second-order model $m_k(x_k + d) $ as follows:
	\begin{equation}
	m_k(x_k + d) = f_k + g_k^{T}d_k + \frac{1}{2}d^{T}B_kd,
	\end{equation}
	where  $f_k$ and $g_k$ are defined as the beginning of this section.
	Next, we will show how to update the quasi-Newton matrix $B_k$.
	
	Let $s_k \stackrel{\Delta}{=} x_{k+1} - x_{k}$ and $y_k \stackrel{\Delta}{=} \frac{1}{b_k}(\sum_{j=1}^{b_k}\nabla f_{i_j}(x_{k+1}) - \nabla f_{i_j}(x_{k})).$ Given a initial matrix $B_0$, provided that $(y_k - B_ks_k)^{T}s_k \neq 0 $, then $B_{k+1}$ can be defined as 
	\begin{equation}
	B_{k+1} \stackrel{\Delta}{=} B_k +  \frac{(y_k - B_ks_k)(y_k - B_ks_k)^{T}}{(y_k - B_ks_k)^{T}s_k}.
	\end{equation}
	Limited memory symmetric rank one (SR1) method stores and uses the $m$ most recently computed pair $\left\lbrace (s_k, y_k)\right\rbrace $ ($m \ll d$). To describe the compact representation of a L-SR1 matrix, we need to define, for $k \geq m$:
	\begin{equation}
	\begin{split}
	S_k & \stackrel{\Delta}{=} \left[ s_{k-m+1}, s_{k-m}, \cdots, s_k\right] \in \R^{d \times m},  \\
	Y_k & \stackrel{\Delta}{=} \left[ y_{k-m+1}, y_{k-m}, \cdots, y_k\right] \in \R^{d \times m}. 
	\end{split}
	\end{equation}
	Moreover, we need the following decomposition of $S_k^{T}Y_k:$
	\begin{equation}
	S_k^{T}Y_k = L_k + D_k + R_k, 
	\end{equation}
	where $L_k$ is strictly lower triangular, $D_k$ is diagonal, and $R_k$ is strictly upper triangular. We assume that all the updates are well-defined, that is $s_k^{T}(y_k - B_ks_k) \neq 0$, otherwise we skip the update. The compact form of L-SR1 \cite{Compact_LSR1} can be written as
	\begin{equation}
	B_{k+1} = B_0 + U_k V_k U_k^{T},
	\end{equation}
	where $U_k \in \R^{d\times (k+1)}$, $V_k \in \R^{(k+1)\times (k+1)}$, and $B_0$ is a diagonal matrix. $U_k$ and $V_k$ are given by
	\begin{equation}
	U_k = Y_k - B_0 S_k \quad \text{and} \quad  V_k = (D_k + L_k + L_k^{T} - S_k^{T}B_0S_k)^{-1}.
	\end{equation}
	Now the quadratic probabilistic model defined by the L-SR1 method is constructed. The trust-region subproblem will be
	\begin{equation}
	\min_{ \left\|d \right\| \leq \delta_k } m_k(x_k + d) \stackrel{\Delta}{=} f_k + g_k^{T}d + \frac{1}{2}d^{T}B_kd.
	\end{equation} 
	Applying this model into STRME, we can obtain  STRME-Lsr1 algorithm. In the same way, we can obtain STORM-Lsr1 algorithm.	 
	With respect to the specific implementation of the subproblem and how to solve the trust-region subproblem efficiently, one can refer to the OBS method in \cite{OBS}.

	In this part, we compare our algorithm STRME-Lsr1 with STORM-Lsr1 and its first-order form. The result is shown in Figure \ref{fig:5}. We test on different choice of matrix $B_0$: (i) $B_0 = \frac{s_0^{T}y_0}{s_0^{T}s_0}$; (ii) $B_0 = \frac{y_0^{T}y_0}{s_0^{T}y_0}$; (iii) $B_0 = \tau_0 \I_d$, where $\tau_0$ is adjusted with the same scale as the initial value  of $\frac{y_0^{T}y_0}{s_0^{T}y_0}$. Unfortunately, the first two choices are not satisfactory. Thus, we use (iii) to define the initial matrix $B_0$. 
	 \vspace{-0.15in}
	\begin{figure}[H]
		\centering
		\includegraphics[width=0.75\textwidth]{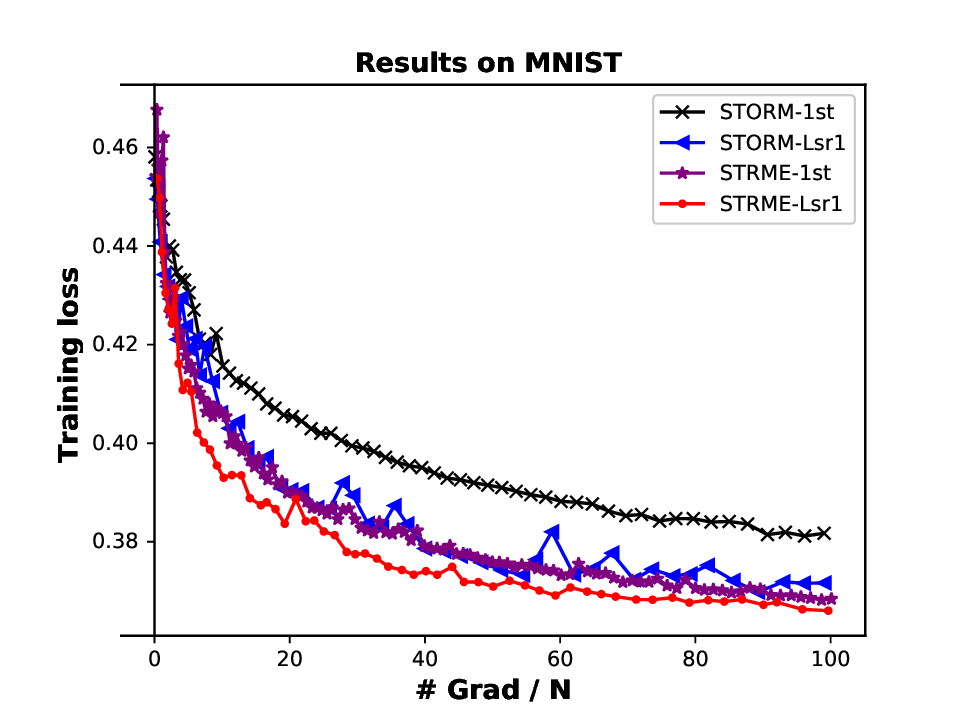}
		\includegraphics[width=0.75\textwidth]{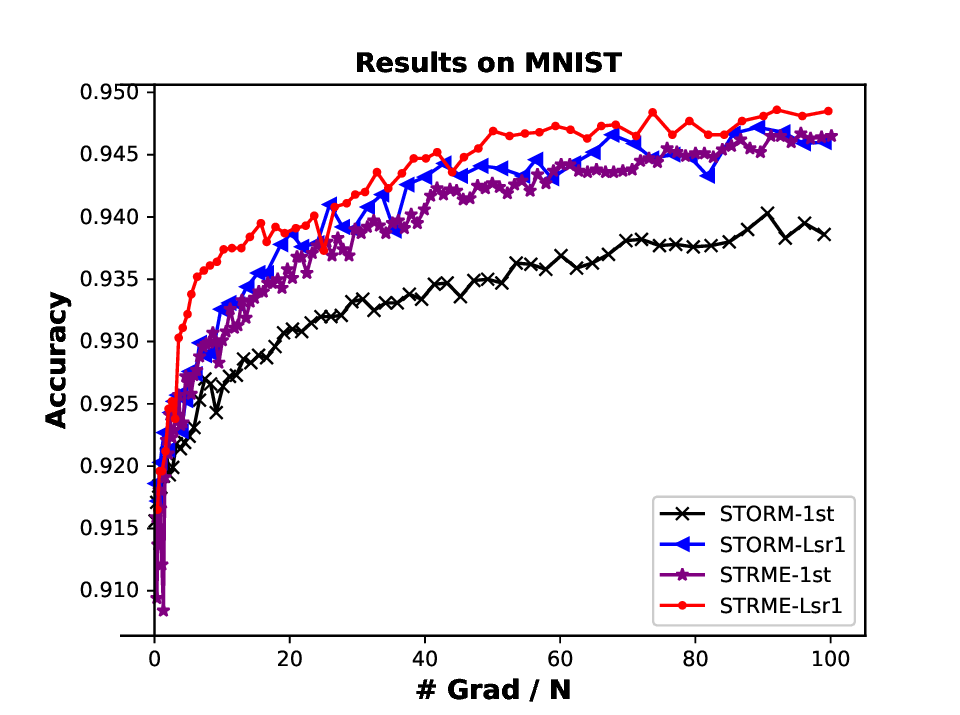}
		\caption{Training DNN on second-order probabilistic model with $ \eta_1=0.1, \gamma=2, m=30, \mu_0 = 0.1, \mu_{\max}=10$ for STRME-Lsr1;  for STORM-Lsr1, $\eta_1=0.1, \eta_2 = 0.001, \gamma=2, m=30, \delta_0=0.1, \delta_{\max} = 1.$ }
		\label{fig:5}
	\end{figure}
	 In our numerical experiement, $B_0 = \tau_0 \I_d$, where $\tau_0 \in \left\lbrace 0.5, 1, 1.5\right\rbrace $, and the best $\tau_0$ is achieved at $\tau_0 = 1$. We have tested the limited memory size $m \in \left\lbrace 10,20,30,40\right\rbrace $. We find that $m=30$ is performed relatively better. The parameters $\mu_0$ is chosen from  $\left\lbrace 0.01, 0.1, 1\right\rbrace $, and the best $\mu_0$ is achieved at $\mu_0 = 0.1$. For $\mu_{\max}$, we test on the range $\left\lbrace 0.1, 1, 10, 100 \right\rbrace $, the best one is $\mu_{\max}=10$, which implies that unlike the first-order STRME, the second-order STRME is not so sensitive to $\mu_{\max}$ and allow some larger values. For STORM-Lsr1, we  test the $\delta_0$ and $\delta_{\max}$ with the same range as $\mu_{0}$ and $\mu_{\max}$, the best choices are achieved at $\delta_0 = 0.1, \delta_{\max}=1$, respectively. From Figure \ref{fig:5}, we can see that our algorithm STRME-Lsr1 performs better than STORM-Lsr1.	 And it is not hard to find that STRME-1st is not worse than STORM-Lsr1.	
	 \vspace{-0.2in}	
	\begin{figure}[H]
		\centering	
		\includegraphics[width=0.75\textwidth]{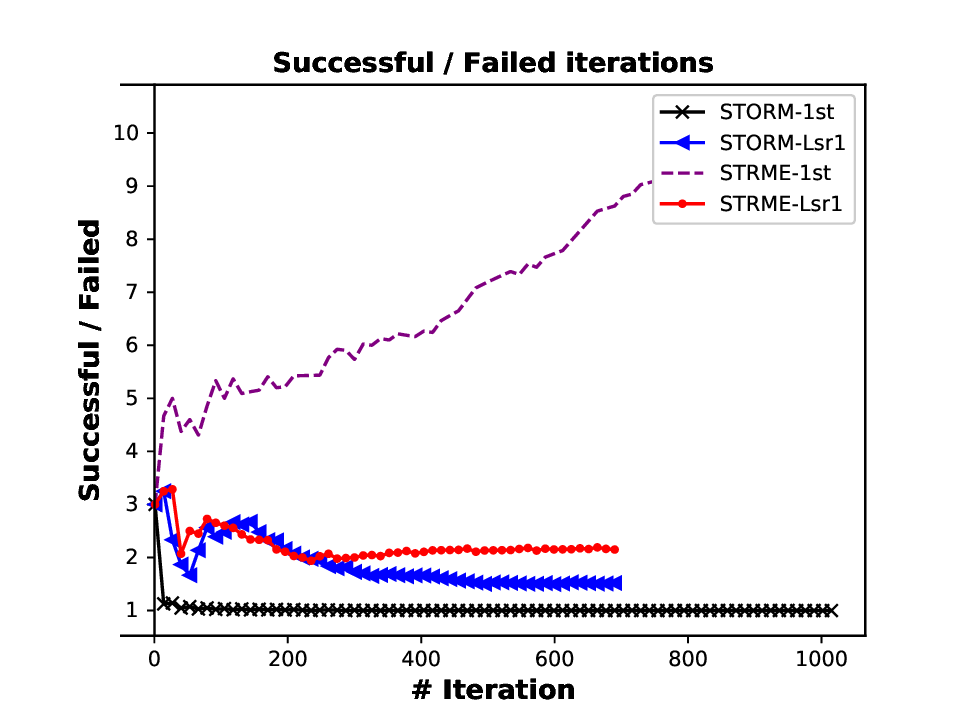}
		\caption{The comparison between STRME-Lsr1 and STORM-Lsr1 on  behavior of successful and failed ratio $\varsigma $ }
		\label{fig:6}
	\end{figure}
	\vspace{-0.2in}		
	\begin{figure}[H]	
		\centering	
		\includegraphics[width=0.75\textwidth]{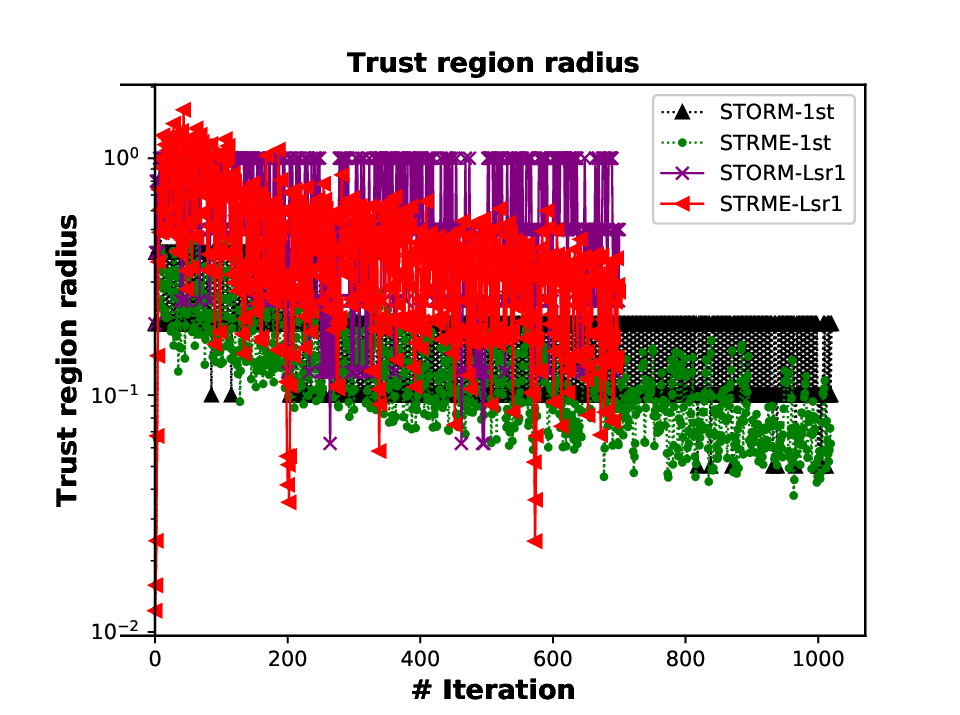}
		\caption{The comparison between STRME-Lsr1 and STORM-Lsr1 on trust-region radius }		
		\label{fig:7}
	\end{figure}
	\vspace{-0.2in}	
				
	Moreover, in Figure \ref{fig:6}, we show the behavior of successful and failed ratio $\varsigma $ which denotes the total number of successful iterations divided by the total number of failed iterations. In this case, we set the maximum number of gradient $SFO_{\max}=100N.$ For the first-order methods, we can see that the ratio $\varsigma $ of STORM-1st is basically stable around 1, and  the ratio of STRME-1st is higher and increasing in the later period. Besides, we can see that in the middle and later period of the algorithm STRME-Lsr1 and STORM-Lsr1, the value of $\varsigma $ is basically stable. However, for STRME-Lsr1, the value of $\varsigma $ in stable condition is still larger than that for STORM-Lsr1. 
	
	Beyond that we also compare the trust-region radius $\delta_k$ of the two algorithms with a long time training in Figure \ref{fig:7}, where $ SFO_{\max} = 100N $. 
	%%%%%%%%%%%%%%%%%%%%%%%%%%%%
	%%%%%%%%%%%%%%这里需要说清楚，分析清楚，因为这里是数值实验的精髓部分。
	We can see that for STRME, whether STRME-1st or STRME-Lsr1, the bandwidth of trust-region radius is narrower than that of STORM. This means the oscillation of the trust-region radius for STRME is less severe in contrast to that of STORM methods. Besides, for STRME-Lsr1, $\delta_k$ is overall declining and smaller than that in STORM-Lsr1. At the same time, we observe that the second-order methods permit larger trust-region radius from the Figure \ref{fig:7}. 
	\vspace{-0.2in}		
	\begin{figure}[H]	
		\centering	
		\includegraphics[width=0.75\textwidth]{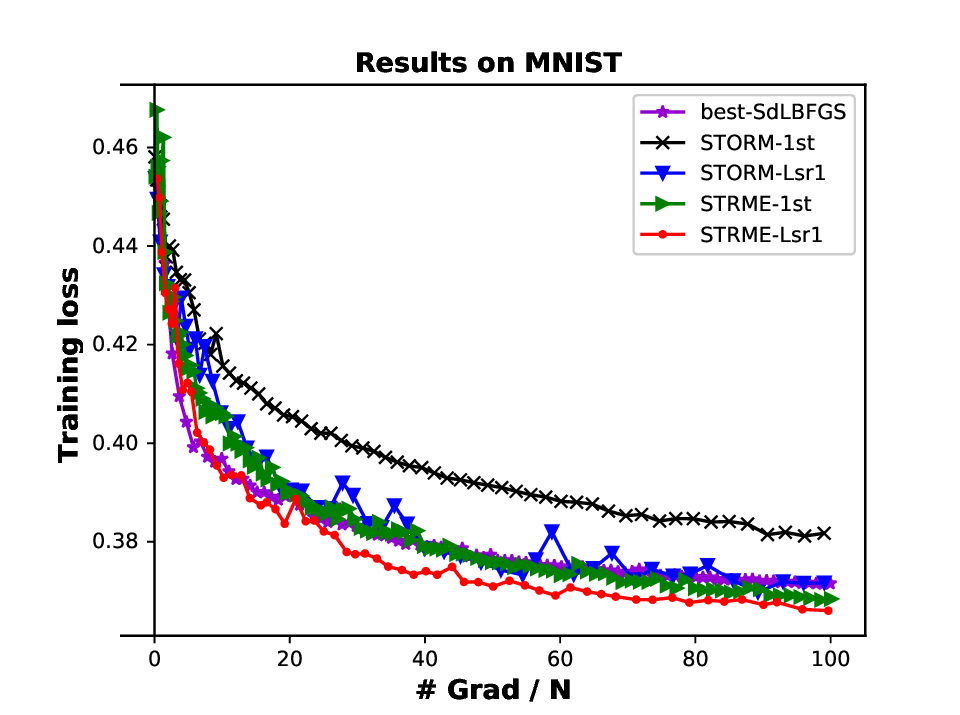}
		
		\includegraphics[width=0.75\textwidth]{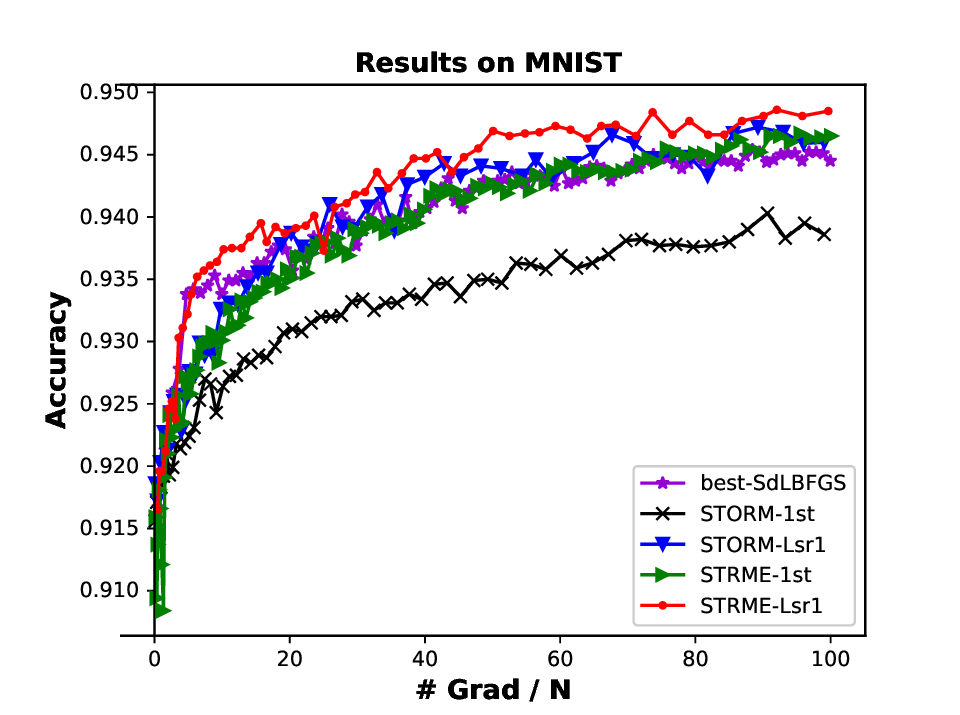}
		\caption{ The Comparison between STRME, SdLBFGS and STORM}
		\label{fig:8}
	\end{figure}
	
	Moreover, we implement the well-known second-order algorithm SdLBFGS \cite{SdLBFGS} to test the performance of our algorithm with a long time running ($ SFO_{\max} = 100N $). In this case, we set batch size $b= d+1$ and step size  $\eta_k = \frac{\beta_0}{(k/10 + 1)}$ where $\beta_0 \in \left\lbrace 0.1, 1, 10 \right\rbrace $ for SdLBFGS. The best tuned step size is obtained at $\beta_0 =10$. In \cite{SdLBFGS}, they set $\eta_k = \frac{10}{k+1}$. By numerical comparison, we find that the result of $\eta_k = \frac{10}{(k/10 + 1)}$ is better than that for the choice $\eta_k = \frac{10}{k + 1}$. The results in Figure \ref{fig:8} illustrate that STRME-Lsr1 performs better than the best tuned SdLBFGS.

	\section{Conclusion}
	\label{sec:5}
	We have presented a stochastic trust-region framework in which the trust-region radius depends on the currently probabilistic model. To verify the effectiveness of the framework, we have proposed a specific algorithm named STRME in which the trust-region radius is linearly associated with the model gradient. We have analyzed the expected number of iterations of STRME for three different cases: non-convex, convex, and strongly convex. We can see that our algorithm enjoys the same complexity properties as the existing schemes. Moreover, our algorithm compares favorably to STORM algorithm and other stochastic algorithms on several testing problems involving the real datasets. Actually, in addition to STRME, there are many other approaches to explore the trust-region radius related with random models. We point out that the work in this paper is limited to the case that the objective function is smooth.  There are some important and latest works for non-smooth problems, for example \cite{RSPG, prox-SVRG-nonconvex,QuickeNing,StSR1}. It is  worthwhile to extend the stochastic trust-region framework to the non-smooth setting. Moreover, the effectiveness of the stochastic trust-region method is relevant to the model. How to construct a more efficient model is an interesting subject for future research.

	\section{Appendix }
	\label{sec:6}

	\subsection*{{\bf A: Proofs of lemmas and theorems in Section 3}}
	\label{sec:6:1}
	{\bf Proof of Lemma \ref{lem1}}
	
	\begin{proof}
		From Assumption \ref{sufficient_reduction}, the trial step $d_k$ will lead a sufficient reduction on $m_k$ such that
		\begin{equation}\label{B_proof_inequ1}
		m_k(x_k) - m_k(x_k + d_k) \geq \frac{\kappa_{fcd}}{2}\left\| g_k\right\| \min\left\lbrace \frac{\left\|g_k \right\| }{\left\|B_k \right\| }, \delta_k\right\rbrace. 
		\end{equation} 
		Since $ \mu_k \leq \frac{1}{\kappa_{bhm}} $, we have $\delta_k = \mu_k\left\|g_k \right\| \leq \frac{\left\|g_k \right\|}{\kappa_{bhm}} \leq \frac{\left\|g_k \right\|}{\left\|B_k \right\| } $ which together with (\ref{B_proof_inequ1}), implies that
		\begin{equation}
		m_k(x_k) - m_k(x_k + d_k) \geq \frac{\kappa_{fcd}}{2}\left\| g_k\right\|^2\mu_k. 
		\end{equation}
		
		Suppose that model $m_k$ is true, we can obtain that
		\begin{equation}
		\begin{split}
		&   f(x_k) - f(x_k + d_k)\\
		= &  f(x_k) - m_k(x_k) + m_k(x_k) - m_k(x_k + d_k)
		+ m_k(x_k+ d_k) - f(x_k + d_k) \\
		\geq &   -2\kappa_{ef}\left\| g_k\right\|^2\mu_k^2 + \frac{\kappa_{fcd}}{2}\left\| g_k\right\|^2\mu_k \\
		= & (-2\kappa_{ef}\mu_k + \frac{\kappa_{fcd}}{2} )\left\| g_k\right\|^2\mu_k.
		\end{split}
		\end{equation}
		Because of the condition that $\mu_k \leq \frac{\kappa_{fcd}}{8\kappa_{ef}}$, we have 
		\begin{equation}
		f(x_k) - f(x_k + d_k) \geq \frac{\kappa_{fcd}}{4}\left\| g_k\right\|^2\mu_k.
		\end{equation}
		Thus, the desired result is proved.
	\end{proof}	
	{\bf Proof of Lemma \ref{lem2}}
	\begin{proof}
		If $d_k$ is accepted, which implies that $\rho_k \geq \eta_1$, then 
		\begin{equation}
		\begin{split}
		f_k^0 -f_k^d & \geq \eta_1(m_k(x_k) - m_k(x_k + d_k)) \\
		& \geq \eta_1 \frac{\kappa_{fcd}}{2}\left\|g_k \right\|\min\left\lbrace \frac{\left\|g_k \right\|}{\left\|B_k \right\| }, \delta_k  \right\rbrace   \\
		& \geq \eta_1 \frac{\kappa_{fcd}}{2}\left\|g_k \right\|^2\min\left\lbrace \frac{1}{\left\|B_k \right\| }, \mu_k  \right\rbrace \\
		& \geq \eta_1 \frac{\kappa_{fcd}}{2}\left\|g_k \right\|^2\mu_k,
		\end{split}
		\end{equation}
		where the last inequality follows from $\mu_k \leq \frac{1}{\kappa_{bhm}}$. 
		
		If the estimates $\left\lbrace f_k^0, f_k^d\right\rbrace $ are tight, the improvement in $f$ can be bounded by
		\begin{equation}
		\begin{split}
		&  f(x_k) - f(x_k+d_k) \\
		= & f(x_k) - f_k^0 + f_k^0 - f_k^d + f_k^d - f(x_k + d_k) \\
		\geq & -2\epsilon_{F}\delta_k^2 + \eta_1 \frac{\kappa_{fcd}}{2}\left\|g_k \right\|^2\mu_k \\
		\geq & (-2\epsilon_{F}\mu_k^2 + \frac{\eta_1\kappa_{fcd}}{2}\mu_k)\left\|g_k \right\|^2 \\
		\geq &  \frac{\eta_1\kappa_{fcd}}{4} \left\|g_k \right\|^2\mu_k.
		\end{split}
		\end{equation}
		Since $\epsilon_F \leq \frac{\eta_1\kappa_{fcd}}{8\mu_{\max}}$, we know that $\epsilon_{F} \leq \frac{\eta_1\kappa_{fcd}}{8\mu_k} $, which deduces the last inequality.
		
	\end{proof}

	{\bf{Proof of Lemma \ref{lem3}}}
	
	\begin{proof}
		Because $\mu_k \leq \frac{1}{\kappa_{bhm}}$, the $Cauchy \, decrease$ condition yields
		\begin{equation}{\label{lem3_ineq_1}}
		m_k(x_k) - m_k(x_k + d_k) \geq \frac{\kappa_{fcd}}{2}\left\| g_k\right\| \min\left\lbrace \frac{\left\|g_k \right\| }{\left\|B_k \right\| }, \delta_k\right\rbrace \geq \frac{\kappa_{fcd}}{2} \left\|g_k \right\|^2 \mu_k.
		\end{equation} 
		Assume that model $m_k$ are true, which means, for all $y \in B(x_k, \delta_k)$, we have
		\begin{equation}{\label{lem3_ineq_2}}
		\begin{split}
		\left\|\nabla f(y) - \nabla m_k(y) \right\|  & \leq  \quad \kappa_{eg}\delta_k, \,\,\,  \text{and} \\
		\left|f(y) - m_k(y) \right|   & \leq \quad  \kappa_{ef} \delta_k^2.	 		
		\end{split}
		\end{equation}		
		And the estimates $f_k^0$ and $f_k^d$ are tight with $\epsilon_{F} \leq \kappa_{ef}$, we have 
		\begin{equation}{\label{lem3_ineq_3}}	
		\left| f_k^0 - f(x_k) \right|   \leq  \quad \kappa_{ef}			
		\delta_k^2, \,\,\,  \text{and} \,\, 	
		\left|f_k^d - f(x_k + d_k) \right|    \leq \quad  \kappa_{ef} \delta_k^2.	 
		\end{equation} 	
		
		The ratio $\rho_k$ can be rewritten as
		\begin{equation}{\label{lem3_equ_4}}
		\begin{split}
		\rho_k  = &  \frac{f_k^0 - f_k^d}{ m_k(x_k) -  m_k(x_k + d_k)} \\ 
		=  &   \frac{f_k^0 - f(x_k)}{ m_k(x_k) -  m_k(x_k + d_k)} + \frac{f(x_k) - m_k(x_k)}{ m_k(x_k) -  m_k(x_k + d_k)} + \frac{m_k(x_k) - m_k(x_k + d_k)}{ m_k(x_k) -  m_k(x_k + d_k)} \\
		& + \frac{m_k(x_k+ d_k) - f(x_k + d_k)}{ m_k(x_k) -  m_k(x_k + d_k)} + \frac{f(x_k + d_k) - f_k^d}{ m_k(x_k) -  m_k(x_k + d_k)}.\\
		\end{split}
		\end{equation}	
		Applying the inequlaities (\ref{lem3_ineq_1}), (\ref{lem3_ineq_2}) and (\ref{lem3_ineq_3}) to the above equality, and then using $\delta_k = \mu_k \left\| g_k \right\|$, we can obtain
		\begin{equation*}
		\left|\rho_k -1 \right| \leq \frac{8\kappa_{ef}\delta_k^2}{\kappa_{fcd}\mu_k \left\| g_k \right\|^2} = \frac{8\kappa_{ef}}{\kappa_{fcd}}\mu_k.
		\end{equation*}	
		Since $\mu_k \leq \frac{\kappa_{fcd}(1-\eta_1)}{8\kappa_{ef}}$, we have $ \frac{8\kappa_{ef}}{\kappa_{fcd}} \mu_k \leq 1-\eta_1$. Thus, we conclude that $\rho_k \geq \eta_1$, which means that the $k$-th iteration is successful.
		
	\end{proof}

	{\bf{Proof of Lemma \ref{lem_assump2(2)}}}	
	\begin{proof}
		First, we show that for all  $k < T_{\epsilon}$, the following inequality holds
		\begin{equation}\label{Lambda_inequ}
		\Lambda_{k+1} \geq \min\left\lbrace \hat{\Lambda}, \min\left\lbrace \mu_{\max}, \gamma \Lambda_k \right\rbrace I_kJ_k + \frac{1}{\gamma}\Lambda_k(1-I_kJ_k)  \right\rbrace.
		\end{equation} 
		If $\Lambda_k > \hat{\Lambda}$, we have $\Lambda_{k+1} \geq \gamma \hat{\Lambda}$ by the update process of the sequence $\Lambda_k$. Hence, $\Lambda_{k+1} \geq \hat{\Lambda}$. Now we assume that $\Lambda_{k} \leq \hat{\Lambda}$, by the definition of $\hat{\Lambda}$, we have 
		\begin{equation}
		\Lambda_{k} \leq \min \left\lbrace \frac{1}{\kappa_{bhm}}, \frac{\kappa_{fcd}(1-\eta_1)}{8\kappa_{ef}} \right\rbrace.
		\end{equation}
		If $I_k = 1$ and $J_k =1$, i.e. model $M_k$ and estimates are all sufficiently accurate,
		from Lemma \ref{lem3}, we know that the iteration $k$ is successful. Thus $X_{k+1} = X_{k} + D_k$ and $\Lambda_{k+1} = \gamma \Lambda_k$. If $I_kJ_k = 0$, whether the iteration $k$ is successful or failed, we all have $\Lambda_{k+1} \geq \frac{1}{\gamma}\Lambda_k$. 
		
		From the above analysis, we conclude that (\ref{Lambda_inequ}) holds. Moreover we have observed that $p = P(W_k =1 |\mathcal{F}_{k-1}^{M,F} ) = P(I_kJ_k =1 | \mathcal{F}_{k-1}^{M,F}) \geq \alpha\beta > \frac{1}{2} $, which implies that Assumption \ref{section2_assump2}(\ref{assump2(2)}) holds.

	\end{proof}

	{\bf{Proof of Theorem \ref{thm_Phi_decrease}}}		
	\begin{proof}
		First, we recall the definition of $\Phi_k$, that is
		\begin{equation}
		\Phi_k = \nu(f(X_k)- f_{\min}) + (1-\nu)\frac{1}{L^2}\Lambda_k\left\|\nabla f(X_k) \right\|^2. 
		\end{equation}
		In the following proof, we consider three separate cases: (i) model $M_k$ is true and estimates $\left\lbrace F_k^0, F_k^d\right\rbrace $ are tight ($I_k=J_k =1$); (ii) model $M_k$ is false and estimates $\left\lbrace F_k^0, F_k^d\right\rbrace $ are tight ($J_k=1, I_k=0$); (iii) estimates $\left\lbrace F_k^0, F_k^d\right\rbrace $ are loose ($J_k=0$). For each of these cases, we will analyze two possible outcomes of the iteration process, that is the iteration $k$ is successful or failed. Based on the above classifications, we rewrite the expected decrease of $\Phi_k$ as
		\begin{equation}
		\E[ \Phi_{k+1} - \Phi_k | \mathcal{F}_{k-1}^{M,F}] = \E[(\I_{\left\lbrace I_kJ_k=1\right\rbrace}  + \I_{\left\lbrace (1-I_k)J_k=1\right\rbrace } + \I_{\left\lbrace (1-J_k)=1\right\rbrace })(\Phi_{k+1} - \Phi_{k}) | \mathcal{F}_{k-1}^{M,F} ].
		\end{equation}
		
		Before presenting the formal proof, we brief describe the key ideas. By choosing a suitable constant $\nu$, we can derive an upper bound on the expected decrease of $\Phi_k$ for each of the cases. When model $M_k$ is true and estimates $\left\lbrace F_k^0, F_k^d\right\rbrace $ are tight, no matter whether the iteration $k$ is successful or not, it will give rise to the decrease of $\Phi_k$ which is in proportion to $\Lambda_k\left\|\nabla f(X_k) \right\|^2$. For the other two cases, due to the model error or inaccurate estimates, $\Phi_k$ may increase. However, the increment of $\Phi_k$ is still in proportion to $\Lambda_k\left\|\nabla f(X_k) \right\|^2$. Therefore, by choosing sufficiently large $\alpha$ and $\beta$,  the expectation of $\Phi_k$ can be guaranteed to decrease.

		\begin{itemize}
			\item[(i).] Model $M_k$ is true and estimates $\left\lbrace F_k^0, F_k^d\right\rbrace $ are tight ($I_k=J_k=1$).
			
			\begin{itemize}
				\item[a.] Iteration $k$ is successful. In this case, we have $X_{k+1} = X_k + D_k$, and $\Lambda_{k+1} = \gamma\Lambda_k.$
				
				Because model $M_k$ is true.
				From Lemma \ref{lem1}, we know that if $$\mu_{\max} \leq \min\left\lbrace \frac{1}{\kappa_{bhm}}, \frac{\kappa_{fcd}}{8\kappa_{ef}} \right\rbrace, $$  we have 
				\begin{equation}\label{inequ1_thm_Phi}
				f(X_{k} + D_k) - f(X_k) \leq - \frac{\kappa_{fcd}}{4}\Lambda_k \left\|G_k \right\|^2. 
				\end{equation}
				Besides, we can easily find the relation between $\left\|G_k \right\| $ and $\left\|\nabla f(X_k) \right\| $, i.e.
				\begin{equation*}
				\begin{split}
				\left\|\nabla f(X_k) - G_k \right\| &  \leq \kappa_{eg}\delta_k = \kappa_{eg} \Lambda_k \left\|G_k \right\|.\\
				\end{split}
				\end{equation*} 
				Using the above inequality, the triangle inequality and the fact that $\Lambda_k \leq \mu_{\max}$, we obtain
				\begin{equation}\label{g_nabla_f}
				\left\|  G_k \right\| \geq \frac{1}{1+\kappa_{eg}\Lambda_{\max}}\left\|\nabla f(X_k) \right\|. 
				\end{equation} 
				Since $f$ is $L$-smooth and iteration $k$ is successful, we have that
				\begin{equation}
				\left\|\nabla f(X_{k+1}) - \nabla f(X_k) \right\| \leq L\Delta_k. 
				\end{equation}
				Using the fact that $ (a+b)^2 \leq 2(a^2 + b^2)$ and by a simple calculation, we get
				\begin{equation}
				\left\|\nabla f(X_{k+1}) \right\|^2 \leq 2(L^2\Lambda_k^2\left\|G_k \right\|^2 + \left\|\nabla f(X_k) \right\|^2).
				\end{equation}
				Particularly, the following holds with $\Lambda_{k+1} = \gamma \Lambda_k$ and $\,\Lambda_k \leq \mu_{\max}$
				\begin{equation}\label{G_nabla_plus}
				\begin{split}
				& \frac{1}{L^2}(\Lambda_{k+1}\left\|\nabla f(X_{k+1}) \right\|^2 - \Lambda_k \left\|\nabla f(X_k) \right\|^2 ) \\
				& \leq 2\gamma\Lambda_k(\mu_{\max}^2\left\|G_k \right\|^2 + \frac{1}{L^2}\left\|\nabla f(X_k) \right\|^2  ).						
				\end{split}
				\end{equation}			
				Applying (\ref{inequ1_thm_Phi}) and (\ref{G_nabla_plus}), we get	
				\begin{equation}\label{inequ2_thm_Phi}
				\begin{split}
				& \Phi_{k+1} - \Phi_k \\
				=  &  \nu(f(X_{k+1} - f(X_k))) + (1-\nu)\frac{1}{L^2}(\Lambda_{k+1}\left\|\nabla f(X_{k+1}) \right\|^2 -\Lambda_{k}\left\|\nabla f(X_{k}) \right\|^2 )  \\
				\leq	& - \frac{\nu\kappa_{fcd}}{4} \Lambda_k \left\|G_k\right\|^2 + 2(1-\nu)\gamma\Lambda_k(\mu_{\max}^2\left\|G_k \right\|^2 + \frac{1}{L^2}\left\|\nabla f(X_k) \right\|^2). \\
				\end{split}
				\end{equation}
				We choose $\nu \in (0,1)$ to satisfy $$\frac{1-\nu}{\nu} \leq \frac{\kappa_{fcd}}{16\gamma \mu_{\max}^2},$$ which implies that
				\begin{equation}
				-\frac{\nu\kappa_{fcd}}{4} \Lambda_k \left\|G_k\right\|^2  + 2(1-\nu)\gamma \mu_{\max}^2\Lambda_k\left\|G_k \right\|^2 \leq -\frac{\nu\kappa_{fcd}}{8}\Lambda_k\left\|G_k \right\|^2.
				\end{equation} 		
				Consequently, (\ref{inequ2_thm_Phi}) can be written as 
				\begin{equation}\label{inequ3_thm_Phi}
				\begin{split}
				\Phi_{k+1} - \Phi_k & \leq -\frac{\nu\kappa_{fcd}}{8}\Lambda_k\left\|G_k \right\|^2 + 2(1-\nu)\gamma\frac{1}{L^2}\Lambda_k \left\|\nabla f(X_k) \right\|^2.
				\end{split}
				\end{equation}
				Applying (\ref{g_nabla_f}) to (\ref{inequ3_thm_Phi}), we have
				\begin{equation}\label{1a_Phi1}
				\begin{split}			
				\Phi_{k+1} - \Phi_k 
				\leq  -\frac{\nu\kappa_{fcd}}{8(1+\kappa_{eg}\mu_{\max})^2}\Lambda_k\left\|\nabla f(X_k) \right\|^2 + 2(1-\nu)\gamma\frac{1}{L^2}\Lambda_k \left\|\nabla f(X_k) \right\|^2.
				\end{split}
				\end{equation}
				Furthermore, we assume that $\nu$ satisfies $$\frac{1-\nu}{\nu} \leq \frac{\kappa_{fcd}L^2}{32\gamma(1+\kappa_{eg}\mu_{\max})^2},$$ which yields 
				\begin{equation}\label{1a_Phi2}
				\begin{split}			
				& -\frac{\nu\kappa_{fcd}}{8(1+\kappa_{eg}\mu_{\max})^2}\Lambda_k\left\|\nabla f(X_k) \right\|^2 + 2(1-\nu)\gamma\frac{1}{L^2}\Lambda_k \left\|\nabla f(X_k) \right\|^2 \\
				\leq &  -\frac{\nu\kappa_{fcd}}{16(1+\kappa_{eg}\mu_{\max})^2}\Lambda_k\left\|\nabla f(X_k) \right\|^2.
				\end{split}
				\end{equation}
				(\ref{1a_Phi1}) and (\ref{1a_Phi2}) show that
				\begin{equation}
				\Phi_{k+1} - \Phi_k \leq -\frac{\nu\kappa_{fcd}}{16(1+\kappa_{eg}\mu_{\max})^2}\Lambda_k\left\|\nabla f(X_k) \right\|^2.
				\end{equation}
				\item[b.] Iteration k is failed. In this case, we know $X_{k+1} = X_k$ and $\Lambda_{k+1} = \frac{1}{\gamma}\Lambda_k$, which deduce that
				\begin{equation}\label{case1_failed}
				\begin{split}		
				\Phi_{k+1} - \Phi_k & = (1-\nu) (\Lambda_{k+1} - \Lambda_k)\frac{1}{L^2}\left\|\nabla f(X_k) \right\|^2  \\ 
				& = -(1-\nu)(1-\frac{1}{\gamma})\frac{1}{L^2}\Lambda_k \left\|\nabla f(X_k) \right\|^2.
				\end{split}
				\end{equation}
			\end{itemize}	
			We can choose a suitable $\nu \in (0,1)$ such that 
			\begin{equation}
			-\frac{\nu\kappa_{fcd}}{16(1+\kappa_{eg}\mu_{\max})^2} \leq -(1-\nu)(1-\frac{1}{\gamma})\frac{1}{L^2}.
			\end{equation}	
			Then,
			\begin{equation}
			\I_{\left\lbrace I_kJ_k=1\right\rbrace }(\Phi_{k+1} - \Phi_{k}) \leq - 	\I_{\left\lbrace I_kJ_k=1\right\rbrace }(1-\nu)(1-\frac{1}{\gamma})\frac{1}{L^2}\Lambda_k \left\|\nabla f(X_k) \right\|^2.
			\end{equation}
			Taking conditional expectation on $\mathcal{F}_{k-1}^{M,F}$, we have
			\begin{equation}\label{inequ_Phi_1}
			\begin{split}
			&\E[\I_{\left\lbrace I_kJ_k=1\right\rbrace }(\Phi_{k+1} - \Phi_{k}) | \mathcal{F}_{k-1}^{M,F}]    \\ 
			\leq  &- P(I_kJ_k=1 | \mathcal{F}_{k-1}^{M,F} )(1-\nu)(1-\frac{1}{\gamma})\frac{1}{L^2}\Lambda_k \left\|\nabla f(X_k) \right\|^2.				
			\end{split}
			\end{equation}

			\item[(ii).] Model $M_k$ is false and estimates $\left\lbrace F_k^0, F_k^d\right\rbrace $ are tight ($J_k =1, I_k=0$).
			
			\begin{itemize}
				\item[a.] Iteration $k$ is successful. Because the iteration is successful, we have $X_{k+1} = X_k + D_k$, and $\Lambda_{k+1} = \gamma\Lambda_k.$	
				In this case, the estimates $\left\lbrace F_k^0, F_k^d\right\rbrace $ are tight. Applying lemma \ref{lem2}, we have
				\begin{equation}\label{inequ4_thm_Phi}
				f(X_{k+1}) - f(X_k) \leq - \frac{\eta_1\kappa_{fcd}}{4}\Lambda_k\left\|G_k \right\|^2,
				\end{equation}
				with $\epsilon_F \leq \frac{\eta_1\kappa_{fcd}}{8\mu_{\max}} $  and $\mu_{\max} \leq \frac{1}{\kappa_{bhm}}$. Due to the fact that $k$ is successful and $f$ is $L$-smooth, so (\ref{G_nabla_plus}) holds. 
				Combining (\ref{G_nabla_plus}) and (\ref{inequ4_thm_Phi}), we have
				\begin{equation}
				\begin{split}
				& \Phi_{{k+1}} - \Phi_{k}\\
				=  & \,\nu (f(X_{k+1}) - f(X_k)) + (1-\nu)\frac{1}{L^2} (\Lambda_{k+1}\left\|\nabla f(X_{k+1}) \right\|^2- \Lambda_k\left\|\nabla f(X_k) \right\|^2 ) \\
				\leq	&  - \frac{\nu\eta_1\kappa_{fcd}}{4}\Lambda_k\left\|G_k \right\|^2 + 2(1-\nu)\gamma\Lambda_k(\mu_{\max}^2\left\|G_k \right\|^2 + \frac{1}{L^2}\left\|\nabla f(X_k) \right\|^2  ).
				\end{split}
				\end{equation}
				Then we choose a suitable $\nu \in (0,1)$ such that $\frac{1-\nu}{\nu} \leq \frac{\eta_1\kappa_{fcd}}{8\gamma\mu_{\max}^2}$, which implies that 
				\begin{equation}
				-\frac{\nu\eta_1\kappa_{fcd}}{4}\Lambda_k\left\|G_k \right\|^2 + 2(1-\nu)\gamma\mu_{\max}^2\Lambda_k\left\|G_k \right\|^2 \leq 0.
				\end{equation}
				Thus, it follows that 
				\begin{equation}
				\Phi_{{k+1}} - \Phi_{k} \leq 2(1-\nu)\gamma\frac{1}{L^2}\Lambda_k\left\|\nabla f(X_k) \right\|^2.
				\end{equation}	
				\item[b.] Iteration $k$ is failed. Here, we have $X_{k+1} = X_k$ and $\Lambda_{k+1} = \frac{1}{\gamma}\Lambda_k$. In this case, (\ref{case1_failed}) holds.
				
			\end{itemize}	
			No matter whether the iteration $k$ is successful or failed, we always have
			\begin{equation}
			\I_{\left\lbrace (1-I_k)J_k = 1\right\rbrace }(\Phi_{k+1} - \Phi_k) \leq 	\I_{\left\lbrace (1-I_k)J_k = 1\right\rbrace }2(1-\nu)\gamma\frac{1}{L^2}\Lambda_k\left\|\nabla f(X_k) \right\|^2.
			\end{equation}
			Taking conditional expectation on the above inequality, we obtain
			\begin{equation}\label{inequ_Phi_2} 
			\begin{split}
			&	\E[\I_{\left\lbrace (1-I_k)J_k=1\right\rbrace }(\Phi_{k+1} - \Phi_k) | \mathcal{F}_{k-1}^{M,F} ]\\ \leq & 2P((1-I_k)J_k=1  | \mathcal{F}_{k-1}^{M,F})(1-\nu)\gamma\frac{1}{L^2}\Lambda_k\left\|\nabla f(X_k) \right\|^2.
			\end{split}
			\end{equation}
			
			\item[(iii).] Estimates  $\left\lbrace F_k^0, F_k^d\right\rbrace $ are loose ($J_k=0$).
			
			\begin{itemize}
				\item[a.] Iteration $k$ is successful. In this case, we have $X_{k+1} = X_k + D_k$, and $\Lambda_{k+1} = \gamma\Lambda_k$. Because iteration $k$ is accepted and the trial step $D_k$ satisfies Assumption \ref{sufficient_reduction}, we have
				\begin{equation}\label{case3_inequ_1}
				\begin{split}
				F_k^0 - F_k^d & \geq \eta_1 (M_k(X_k) - M_k(X_{k+1})) \\
				& \geq \eta_1 \kappa_{fcd}\left\|G_k \right\|\min\left\lbrace \frac{\left\|G_k \right\| }{\left\| B_k\right\| }, \Delta_k \right\rbrace \\
				& \geq \eta_1 \kappa_{fcd}\Lambda_k\left\|G_k \right\|^2, 
				\end{split}
				\end{equation} 
				where the last inequality reuses the fact that $\Delta_k = \Lambda_k \left\|G_k \right\| $ and $\mu_{\max} \leq \frac{1}{\kappa_{bhm}}$. Applying (\ref{case3_inequ_1}), a successful iteration deduces the following bound for the increment of $f$
				\begin{equation}
				\begin{split}
				f(X_{k+1}) - f(X_k) & = f(X_{k+1}) - F_k^d + F_k^d - F_k^0 + F_k^0 - f(X_k) \\
				& \leq \left|f(X_{k+1}) - F_k^d \right| + F_k^d - F_k^0 + \left|F_k^0 - f(X_k)\right|\\
				& \leq - \eta_1 \kappa_{fcd}\Lambda_k\left\|G_k \right\|^2 + \left|f(X_{k+1}) - F_k^d \right| + \left|F_k^0 - f(X_k)\right|. 
				\end{split}
				\end{equation}
				Then, the upper bound for $\Phi_{k+1} - \Phi_k$ can be obtained as follows
				\begin{equation}\label{case3_inequ_2}
				\begin{split}
				\Phi_{k+1} - \Phi_k =  & \,\, \nu (f(X_{k+1}) - f(X_k)) + (1-\nu)\frac{1}{L^2} (\Lambda_{k+1}\left\|\nabla f(X_{k+1}) \right\|^2- \Lambda_k\left\|\nabla f(X_k) \right\|^2 ) \\
				\leq	& \,\, \nu (f(X_{k+1}) - f(X_k)) + 2(1-\nu)\gamma\Lambda_k(\mu_{\max}^2\left\|G_k \right\|^2 + \frac{1}{L^2}\left\|\nabla f(X_k) \right\|^2 ) \\
				\leq	&  \,\,\nu(- \eta_1 \kappa_{fcd}\Lambda_k\left\|G_k \right\|^2 + \left|f(X_{k+1}) - F_k^d \right| + \left|F_k^0 - f(X_k)\right| ) \\
				& + 2(1-\nu)\gamma\Lambda_k(\mu_{\max}^2\left\|G_k \right\|^2 + \frac{1}{L^2}\left\|\nabla f(X_k) \right\|^2 ),
				\end{split}
				\end{equation}
				where the first inequality uses (\ref{G_nabla_plus}), which is still true in this setting. 
				We can choose $\nu \in (0,1)$ to satisfy $\frac{1-\nu}{\nu} \leq \frac{\eta_1\kappa_{fcd}}{4\gamma\mu_{\max}^2}$ so that 
				\begin{equation}
				- \nu\eta_1 \kappa_{fcd}\Lambda_k\left\|G_k \right\|^2 + 2(1-\nu)\gamma\Lambda_k\mu_{\max}^2\left\|G_k \right\|^2 \leq - \frac{1}{2}\nu\eta_1 \kappa_{fcd}\Lambda_k\left\|G_k \right\|^2.
				\end{equation}
				Then (\ref{case3_inequ_2}) can be rewritten as
				\begin{equation}\label{inequ_4_thm_Phi}
				\begin{split}		
				\Phi_{k+1} - \Phi_k	 \leq &   -\frac{1}{2}\nu\eta_1 \kappa_{fcd}\Lambda_k\left\|G_k \right\|^2 + 2(1-\nu)\gamma\frac{1}{L^2}\Lambda_k\left\|\nabla f(X_k) \right\|^2  \\
				& +  \nu(\left|f(X_{k+1}) - F_k^d \right| + \left|F_k^0 - f(X_k)\right|).						
				\end{split}
				\end{equation}
				
				\item[b.] Iteration $k$ is failed. In this case, we have $X_{k+1} = X_k$ and $\Lambda_{k+1} = \frac{1}{\gamma}\Lambda_k$. Then (\ref{case1_failed}) holds.
				
			\end{itemize}
			Compared to (\ref{case1_failed}), we see that (\ref{inequ_4_thm_Phi}) dominates the upper bound of $\Phi_{k+1} - \Phi_k$. Taking conditional expectations on  inequality (\ref{inequ_4_thm_Phi}) and applying Assumption \ref{section2_assump1}(\ref{assump1(3)}), we have
			\begin{equation}
			\begin{split}
			& \E[\I_{\left\lbrace (1-J_k)=1\right\rbrace }(\Phi_{k+1} - \Phi_{k}) | \mathcal{F}_{k-1}^{M,F}]  \\
		\leq	&  P(1-J_k=1| \mathcal{F}_{k-1}^{M,F})[-\frac{1}{2}\nu\eta_1 \kappa_{fcd}\Lambda_k\left\|G_k \right\|^2 + \nu(2\kappa_{f}\Delta_k^2) ] \\ 
			& + 2P(1-J_k=1| \mathcal{F}_{k-1}^{M,F})(1-\nu)\gamma\frac{1}{L^2}\Lambda_k\left\|\nabla f(X_k) \right\|^2 \\
			\leq  & P(1-J_k=1| \mathcal{F}_{k-1}^{M,F})\nu(-\frac{1}{2}\eta_1 \kappa_{fcd} + 2\kappa_{f}\Lambda_k )\Lambda_k\left\|G_k\right\|^2 \\
			& +2P(1-J_k=1| \mathcal{F}_{k-1}^{M,F})(1-\nu)\gamma\frac{1}{L^2}\Lambda_k\left\|\nabla f(X_k) \right\|^2. \\
			\end{split}
			\end{equation}
			From relations $\kappa_f \leq 2\eta_1\kappa_{ef}$  and  $\Lambda_k \leq \mu_{\max} \leq \min\left\lbrace \frac{1}{\kappa_{bhm}}, \frac{\kappa_{fcd}}{8\kappa_{ef}} \right\rbrace$, it follows that 
			\begin{equation}
			-\frac{1}{2}\eta_1 \kappa_{fcd} + 2\kappa_{f}\Lambda_k \leq -\frac{1}{2}\eta_1 \kappa_{fcd} + 2\kappa_{f}\frac{\kappa_{fcd}}{8\kappa_{ef}}\leq 0.
			\end{equation}
			Thus
			\begin{equation}\label{inequ_Phi_3} 
			\begin{split}		
			& \E[\I_{\left\lbrace (1-J_k)=1\right\rbrace }(\Phi_{k+1} - \Phi_{k}) | \mathcal{F}_{k-1}^{M,F}] \\
			\leq & 2P(1-J_k=1 | \mathcal{F}_{k-1}^{M,F})(1-\nu)\gamma\frac{1}{L^2}\Lambda_k\left\|\nabla f(X_k) \right\|^2.
			\end{split}
			\end{equation}
		\end{itemize}	
		
		Now combining (\ref{inequ_Phi_1}), (\ref{inequ_Phi_2}) and (\ref{inequ_Phi_3}), we can show that 
		\begin{equation}\label{all_inequ}
		\begin{split}
		\E[\Phi_{{k+1}} - \Phi_k | \mathcal{F}_{k-1}^{M,F}] =  & \E[ (\I_{\left\lbrace I_kJ_k=1\right\rbrace } +\I_{\left\lbrace (1-I_k)J_k=1\right\rbrace } + \I_{\left\lbrace (1-J_k)=1\right\rbrace })(\Phi_{k+1} - \Phi_{k}) |\mathcal{F}_{k-1}^{M,F} ] \\
		\leq	&  - P( I_kJ_k=1 | \mathcal{F}_{k-1}^{M,F})(1-\nu)(1-\frac{1}{\gamma})\frac{1}{L^2}\Lambda_k \left\|\nabla f(X_k) \right\|^2 \\
		& +  2P( (1-I_k)J_k=1| \mathcal{F}_{k-1}^{M,F})(1-\nu)\gamma\frac{1}{L^2}\Lambda_k\left\|\nabla f(X_k) \right\|^2 \\
		& +  2P( 1-J_k=1 | \mathcal{F}_{k-1}^{M,F})(1-\nu)\gamma\frac{1}{L^2}\Lambda_k\left\|\nabla f(X_k) \right\|^2 \\
		\leq	&  - P( I_kJ_k=1 | \mathcal{F}_{k-1}^{M,F})(1-\nu)(1-\frac{1}{\gamma})\frac{1}{L^2}\Lambda_k \left\|\nabla f(X_k) \right\|^2 \\
		& + 2P(\left\lbrace (1-I_k)J_k=1\right\rbrace  \bigcup \left\lbrace 1-J_k=1\right\rbrace  | \mathcal{F}_{k-1}^{M,F})(1-\nu)\gamma\frac{1}{L^2}\Lambda_k\left\|\nabla f(X_k) \right\|^2 \\
		\leq 	& - \alpha\beta(1-\nu)(1-\frac{1}{\gamma})\frac{1}{L^2}\Lambda_k \left\|\nabla f(X_k) \right\|^2 + 2(1-\alpha\beta)(1-\nu)\gamma\frac{1}{L^2}\Lambda_k \left\|\nabla f(X_k) \right\|^2 \\
		\leq 	& - \frac{1}{2}\alpha\beta(1-\nu)(1-\frac{1}{\gamma})\frac{1}{L^2}\Lambda_k \left\|\nabla f(X_k) \right\|^2.
		\end{split}
		\end{equation}
		By a simple caculation, we have $ P(\I_{\left\lbrace I_kJ_k=1\right\rbrace } | \mathcal{F}_k^{M,F}) \geq \alpha\beta$. The sixth line of above inequality uses the fact that event $ \left\lbrace (1-I_k)J_k = 1\right\rbrace $ and $\left\lbrace (1-J_k) = 1 \right\rbrace $ are disjoint, which implies that $$ \left\lbrace (1-I_k)J_k = 1\right\rbrace \bigcup \left\lbrace (1-J_k) = 1 \right\rbrace   = \left\lbrace 1-I_kJ_k = 1\right\rbrace. $$ Thus, we have $P(\left\lbrace (1-I_k)J_k = 1\right\rbrace  \bigcup \left\lbrace (1-J_k) = 1\right\rbrace  | \mathcal{F}_{k-1}^{M,F}) = P( 1-I_kJ_k = 1 | \mathcal{F}_{k-1}^{M,F}) \leq 1-\alpha\beta$. Choosing suitable $\alpha$ and $\beta$ such that
		\begin{equation}
		- \alpha\beta(1-\frac{1}{\gamma})  +2(1-\alpha\beta)\gamma\leq  -\frac{1}{2}\alpha\beta(1-\frac{1}{\gamma}),
		\end{equation}
		which implies that
		\begin{equation}
		\alpha\beta \geq \frac{4\gamma^2}{4\gamma^2+ (\gamma-1) } < 1,
		\end{equation}
		then we get the last inequality of (\ref{all_inequ}).
		
		Finally, the proof is completed.	
		
	\end{proof}	
	
	{\bf{Proof of Theorem \ref{thm_convex}}}
	\begin{proof}
		In this case, $f$ satisfies Assumption \ref{f_convex}. Let us consider the stochastic process $\left\lbrace \Lambda_k, \Psi_k \right\rbrace $ with 
		\begin{equation}
		\Psi_{k} = \frac{1}{\nu\epsilon} - \frac{1}{\Phi_k}.
		\end{equation}		
		The convexity of $f$ implies that
		\begin{equation}
		f(x) - f(y) \geq \nabla f(y)^{T}(x-y).
		\end{equation}
		Let $x = x^{\ast}$, $y = X_k$, it follows from the above inequality that 
		\begin{equation}\label{f_convex_inequ1}
		\begin{split}
		f(X_k) - f(x^{\ast}) \leq \nabla f(X_k)^{T}(X_k - x^{\ast}) \leq \left\|\nabla f(X_k) \right\|\left\|X_k - x^{\ast} \right\|.  
		\end{split}
		\end{equation}	
		Because $f$ is $L$-smooth, we have $\left\|\nabla f(X_k) - \nabla f(x^{\ast}) \right\| \leq L \left\|X_k - x^{\ast} \right\|  $. Due to Assumption \ref{f_convex}, we know the level set $\mathcal{L} $ is bounded, then
		\begin{equation}\label{f_convex_inequ2}
		\left\|\nabla f(X_k) \right\|\leq LD.
		\end{equation}
		Combining (\ref{f_convex_inequ1}) and (\ref{f_convex_inequ2}), we have  
		\begin{equation}
		\begin{split}
		\Phi_k & = \nu(f(X_k) - f^{\ast}) + (1-\nu)\frac{1}{L^2}\Lambda_k \left\|\nabla f(X_k) \right\|^2 \\
		& \leq \nu\left\|\nabla f(X_k) \right\|\left\|X_k - x^{\ast} \right\| + (1-\nu)\frac{1}{L^2}\mu_{\max}\left\|\nabla f(X_k) \right\|^2   \\ 
		& \leq (\nu  + (1-\nu)\frac{1}{L}\mu_{\max} )D\left\|\nabla f(X_k) \right\|. 
		\end{split}
		\end{equation}
		From the above inequality and the result of Theorem \ref{thm_Phi_decrease}, we have
		\begin{equation}\label{f_convex_Phi}
		\begin{split}
		\E[ \Phi_{k+1} - \Phi_k | \mathcal{F}_{k-1}^{M,F}]  & \leq - \frac{1}{2}\alpha\beta(1-\nu)(1-\frac{1}{\gamma})\frac{1}{L^2}\Lambda_k \left\|\nabla f(X_k) \right\|^2  \\
		& \leq - \frac{\alpha\beta(1-\nu)(1-\frac{1}{\gamma})}{2(\nu L  + (1-\nu)\mu_{\max} )^2D^2}\Lambda_k \Phi_k^2.
		\end{split}
		\end{equation}
		(\ref{f_convex_Phi}) implies that $\E[\Phi_{k+1} | \mathcal{F}_k^{M,F}] \leq \Phi_k$.
		Recalling the definition of $\Psi_k$, for all $k < T_{\epsilon}$, we have
		\begin{equation}\label{f_convex_inequ3}
		\begin{split}
		\E[\Psi_{k+1} - \Psi_k | \mathcal{F}_{k-1}^{M,F}] &  = \E[ \frac{1}{\Phi_{k}} - \frac{1}{\Phi_{k+1}} | \mathcal{F}_{k-1}^{M,F}] \\
		& \leq \frac{1}{\Phi_{k}} -  \frac{1}{\E[\Phi_{k+1} | \mathcal{F}_{k-1}^{M,F}]} = \frac{\E[\Phi_{k+1} |\mathcal{F}_{k-1}^{M,F}] -\Phi_{k} }{\Phi_{k}\E[\Phi_{k+1} | \mathcal{F}_{k-1}^{M,F}]} \\
		& \leq  - \frac{\alpha\beta(1-\nu)(1-\frac{1}{\gamma})}{2(\nu L + (1-\nu)\mu_{\max} )^2D^2} \Lambda_k\frac{\Phi_k^2}{\Phi_{k}\E[\Phi_{k+1} | \mathcal{F}_{k-1}^{M,F}]} \\
		& \leq  - \frac{\alpha\beta(1-\nu)(1-\frac{1}{\gamma})}{2(\nu L  + (1-\nu)\mu_{\max} )^2D^2} \Lambda_k. 
		\end{split}
		\end{equation}
		The first inequality of (\ref{f_convex_inequ3}) follows from Jensen's inequality which will be given in Lemma \ref{jensen_inequ}.  The second inequality uses (\ref{f_convex_Phi}). The last inequality is due to the fact that $\E[\Phi_{k+1} | \mathcal{F}_k^{M,F}] \leq \Phi_k$. Here, we define an non-decreasing function $h(\cdot)$ as follows
		\begin{equation}
		h(\Lambda_k) =  C_1\Lambda_k, 
		\end{equation}
		where $C_1 = \frac{\alpha\beta(1-\nu)(1-\frac{1}{\gamma})}{2(\nu L  + (1-\nu)\mu_{\max} )^2D^2}$. Then we know that Assumption \ref{section2_assump2}(\ref{assump2(3)}) holds.
		
		From Lemma \ref{lem_assump2(2)}, we can easily obtain that if $\alpha\beta > \frac{1}{2}$ and $\hat{\Lambda}$ is defined as (\ref{hat_Lambda}), Assumption \ref{section2_assump2}(\ref{assump2(2)}) satisfies. Then we have Assumption \ref{section2_assump2} holds. Thus, the conclusion of Theorem \ref{thm1} is true in this case.
		
		Finally, substituting the expression of $\Psi_0$, $\hat{\Lambda}$ and $h$ into Theorem \ref{thm1}, we have
		\begin{equation}
		\E[ T_{\epsilon}] \leq \frac{\alpha\beta}{2\alpha\beta-1}(\frac{M}{\epsilon} + \mathcal{O}(1)),
		\end{equation}
		where $ M = \frac{2(\nu L  + (1-\nu)\mu_{\max} )^2D^2}{\alpha\beta\nu(1-\nu)(1-\frac{1}{\gamma})\hat{\Lambda}}.$  Here, we simplify the constant term as $\mathcal{O}(1)$.
		
		Now, this completes the proof of Theorem \ref{thm_convex}.
	\end{proof}
	
	{\bf{Proof of Theorem \ref{thm_stronglyconvex}}}
	\begin{proof}
		In this setting, $f$ is strongly convex with $\sigma > 0$. We will consider the measure $\Psi_k $ as follows 
		\begin{equation}
		\Psi_k = \log(\Phi_k) + \log(\frac{1}{\nu\epsilon}),
		\end{equation}
		to analyze the theoretical complexity.
		Due to the strongly convexity, we have 
		$$f(X_k) - f^{\ast} \leq \frac{1}{2\sigma}\left\|\nabla f(X_k) \right\|^2. $$
		Then,
		\begin{equation}\label{f_stronglyconvex_inequ1}
		\begin{split}
		\Phi_k & = \nu(f(X_k) - f^{\ast}) + (1-\nu)\frac{1}{L^2}\Lambda_k \left\|\nabla f(X_k) \right\|^2 \\
		& \leq  \nu \frac{\left\|\nabla f(X_k) \right\|^2 }{2\sigma} + (1-\nu)\frac{1}{L^2}\mu_{\max}\left\|\nabla f(X_k) \right\|^2 \\
		& \leq (\frac{\nu}{2\sigma}  + (1-\nu)\frac{1}{L^2}\mu_{\max} )\left\|\nabla f(X_k) \right\|^2. 
		\end{split}
		\end{equation}
		It follows from (\ref{f_stronglyconvex_inequ1}) and Theorem \ref{thm_Phi_decrease} that
		\begin{equation}
		\begin{split}
		\E[ \Phi_{k+1} - \Phi_k | \mathcal{F}_{k-1}^{M,F}]  & \leq - \frac{1}{2}\alpha\beta(1-\nu)(1-\frac{1}{\gamma})\Lambda_k \left\|\nabla f(X_k) \right\|^2  \\
		& \leq - \frac{\alpha\beta(1-\nu)(1-\frac{1}{\gamma})}{(\frac{\nu}{2\sigma}  + (1-\nu)\frac{1}{L^2}\mu_{\max} )}\Lambda_k \Phi_k.
		\end{split}
		\end{equation}
		The above inquality implies
		\begin{equation}
		\E[ \Phi_{k+1} | \mathcal{F}_{k-1}^{M,F}] \leq (1- \frac{\alpha\beta(1-\nu)(1-\frac{1}{\gamma})}{(\frac{\nu}{2\sigma}  + (1-\nu)\frac{1}{L^2}\mu_{\max} )}\Lambda_k) \Phi_k. 
		\end{equation}
		Recalling the definition $\Psi_{k} = \log(\Phi_k) + \log(\frac{1}{\nu\epsilon})$, we have
		
		\begin{equation}
		\begin{split}
		\E[\Psi_{k+1} - \Psi_{k} | \mathcal{F}_{k-1}^{M,F}] & = \E[ \log (\Phi_{k+1}) - \log(\Phi_k) | \mathcal{F}_{k-1}^{M,F}] \\
		&  \leq \log(\E[\Phi_{k+1} | \mathcal{F}_{k-1}^{M,F} ]) - \log(\Phi_k) \\
		& \leq \log(1- \frac{\alpha\beta(1-\nu)(1-\frac{1}{\gamma})}{(\frac{\nu}{2\sigma}  + (1-\nu)\frac{1}{L^2}\mu_{\max} )}\Lambda_k) \\
		& \leq - \frac{\alpha\beta(1-\nu)(1-\frac{1}{\gamma})}{(\frac{\nu}{2\sigma}  + (1-\nu)\frac{1}{L^2}\mu_{\max} )}\Lambda_k.
		\end{split} 
		\end{equation}
		We can define 
		\begin{equation}
		h(\Lambda_k) = C_2\Lambda_k, 
		\end{equation}
		where $C_2 = \frac{\alpha\beta(1-\nu)(1-\frac{1}{\gamma})}{(\frac{\nu}{2\sigma}  + (1-\nu)\frac{1}{L^2}\mu_{\max} )}$.
		
		From Lemma \ref{lem_assump2(2)}, we can easily see that Assumption \ref{section2_assump2}(\ref{assump2(2)}) satisfies if $\alpha\beta > \frac{1}{2}$ and $\hat{\Lambda}$ is defined as (\ref{hat_Lambda}). Thus Assumption \ref{section2_assump2} holds. So the conclusion of Theorem \ref{thm1} is true in strongly convex setting.
		
		By substituting the expression of $\Psi_0$, $\hat{\Lambda}$ and $h$ into Theorem \ref{thm1}, we have
		\begin{equation}
		\E[ T_{\epsilon}] \leq \frac{\alpha\beta}{2\alpha\beta-1}(M\log(\frac{1}{\epsilon}) + \mathcal{O}(1)),
		\end{equation}
		where $ M = \frac{(\frac{\nu}{2\sigma}  + (1-\nu)\frac{1}{L^2}\mu_{\max} )}{\alpha\beta(1-\nu)(1-\frac{1}{\gamma})\hat{\Lambda}}.$
		
		Now the proof is finished.
		
	\end{proof}
	\subsection*{{\bf B: Related lemmas and algorithms for second-order STRME }} \label{sec:6:2}
	\begin{lem}[Chebyshev Inequality\cite{durrett2019probability}]
		%\label{Chebyshev_inequ}
		If $X$ is a random variable with mean $\E[X]$  and variance $Var(X)$, then
		\begin{equation}
		P(| X - \E[X] | \geq v) \leq \frac{Var(X)}{v^2}, \,\,\forall  v > 0.
		\end{equation}
	\end{lem}
Based on the exercise 4.1.2 in \cite{durrett2019probability}, we prove the Chebyshev Inequality with conditional expectation.
	
\begin{lem}\label{Chebyshev_inequ}
If $X$ is a random variable	given the $\sigma$-field $\mathcal{F}$, then %with conditional expectation $\E[X | \mathcal{F}]$ and conditional variance $Var[ X \,|\, \mathcal{F} ]$ 
	\begin{equation}
	P(| X - \E[X | \mathcal{F}] | \geq v \,\, | \,\, \mathcal{F}) \leq \frac{Var[ X \,|\, \mathcal{F} ]}{v^2}, \,\,\forall  v > 0.
	\end{equation}
	\begin{proof}
		Let $A = \left\lbrace  X \,|\, | X - \E[X] | \geq v \right\rbrace $, for $\forall v > 0$, then
		\begin{equation}
		\begin{split}
			Var(X | \mathcal{F}) & = \sum_{s}P(s | \mathcal{F})  | X(s) - \E[X] |^2	\\
			& =  \sum_{s \in A}P(s | \mathcal{F})  | X(s) - \E[X] |^2 +\sum_{s \notin A}P(s | \mathcal{F})  | X(s) - \E[X] |^2 \\
			& \geq \sum_{s \in A}P(s | \mathcal{F})  | X(s) - \E[X] |^2\\
			& \geq  \sum_{s \in A}P(s | \mathcal{F}) v^2 \\
			& \geq v^2 \sum_{s \in A}P(s | \mathcal{F}) \\
			& = v^2 P(A | \mathcal{F}) \\
			& = v^2 P(| X - \E[X | \mathcal{F}] | \geq v \,\, | \,\, \mathcal{F}). 
		\end{split}		
		\end{equation}
	Thus, the proof is finished.	
	\end{proof}
\end{lem}

	\begin{lem}[Jensen Inequality\cite{Probability_Bertsekas}]\label{jensen_inequ}		
		Assume that $f$ is continuous and convex. If X is a random variable, then
		\begin{equation}
		\E[f(X)] \geq f(\E[X] ).
		\end{equation}
	\end{lem}
	\begin{rem}
		If $f$ is concave in Lemma \ref{jensen_inequ}, then we get the opposite result, i.e. $\E[f(X)] \leq f(\E[X])$.
	\end{rem}
	\begin{algorithm}[H]
		\caption{STRME with Dogleg for the logistic loss problem}\label{alg:3}
		\begin{algorithmic}[1]			
			\STATE \textbf{Initialization}: initial point $x_0$, $\gamma > 1$, $\eta_1\in(0,1)$, $\mu_0$,  $k = 0$, $SFO_{\max}$, $\epsilon=10^{-8}$, TotalSFO=0, k=0
			\WHILE{ TotalSFO $\leq$ $SFO_{\max}$ }					
			\STATE Compute $g_k = \frac{1}{b_k}\sum_{i\in O_k} \nabla f_i(x_k)$, $B_k = \frac{1}{b_k}\sum_{i\in O_k} \nabla^2 f_i(x_k)$, where the mini-batch set $O_k$ is randomly chosen
			\IF {$\left\|g_k \right\| \leq  \epsilon $}
			\STATE return to step 3 until $g_k >\epsilon$
			\ENDIF	
			\STATE TotalSFO = TotalSFO + $b_k$
			\STATE Compute $\delta_{k} = \mu_{k}\left\|g_{k} \right\| $ 
			\STATE Compute the Cauchy point $d^{u} = - \frac{g_k^{T}g_k}{g_k^{T}B_kg_k}g_k $
			\IF{ $\left\|d^{u} \right\| \geq \delta_k $ } 
			\STATE $d_k = d_k^{u}$ 
			\ELSE 
			\STATE Compute the Newton step $d^{B} = - B_k^{-1}g_k$ 
			\IF{ $\left\| d^{B}\right\| \leq \delta_k $}
			\STATE $d_k = d^{B}$
			\ELSE
			\STATE Compute $t_b$ to satisfy $\left\| d^u + t(d^{B} - d^{u})\right\| = \delta_k$
			\STATE $d_k = d^{u} + t_b (d^{B} - d^{u})$
			\ENDIF            
			\ENDIF
			\STATE Compute Pred = $- (g_k^{T}d_k + \frac{1}{2}d_kB_kd_k)$ 
			\STATE Obtain estimates $f_k^0$ and $f_k^d$ of $f(x_k)$ and $f(x_k + d_k)$
			
			\STATE Compute $\rho_k = \frac{f_k^0 - f_k^d}{Pred}$
			\IF {$\rho_k \geq \eta_1$}
			\STATE  $x_{k+1} = x_k + d_k$
			\STATE $\mu_{k+1} = \min(\gamma\mu_k, \mu_{\max})$
			\ELSE 
			\STATE $x_{k+1} = x_k$	
			\STATE $\mu_{k+1} = \mu_k/\gamma$ 
			\ENDIF		
			
			\STATE Set $k := k+1$
			\ENDWHILE
		\end{algorithmic}
	\end{algorithm}	
	\begin{algorithm}[H]
		\caption{STRME with L-SR1 for the DNN problem}\label{alg:4}
		\begin{algorithmic}[1]		
			\STATE \textbf{Initialization}: initial point $x_0$, $\gamma > 1$, $\eta_1\in(0,1)$, $\mu_0$, $\mu_{\max} $, $SFO_{\max}$, $\epsilon=10^{-6}$, TotalSFO = 0, $t_0$, $b_0, b_{\max}$, $\epsilon=10^{-8}$; Set k=0
			\WHILE{ TotalSFO $\leq$ $SFO_{\max}$ }					
			\STATE Compute $g_k = \frac{1}{b_k}\sum_{i\in O_k} \nabla f_i(x_k)$, where the mini-batch set $O_k$ is randomly chosen without replacement and $b_k = \min \left\lbrace b_{\max}, \max\left\lbrace t_0 k + b_0, \frac{1}{\delta^2} \right\rbrace \right\rbrace $
			\IF {$\left\|g_k \right\| \leq  \epsilon $}
			\STATE return to step 2 until $g_k >\epsilon$
			\ENDIF
			\STATE TotalSFO = TotalSFO + $b_k$
			\STATE Compute $\delta_{k} = \mu_{k}\left\|g_{k} \right\| $ 
			\IF{ len(S)==0}
			\STATE $s_k = - \frac{\delta_k}{\left\|g_k \right\| }g_k$, and $ B_k = B_0$
			\ELSE
			\STATE Update $B_{k+1}$ as in section \ref{sec:4:3:2}
			\STATE $s_k = \arg\min_{\left\| s\right\| \leq \delta_k} m_k(x_k + s)$ 
			\STATE Compute $B_ks_k$
			\ENDIF
			\STATE Compute Pred $= - (g_k^{T}s_k + \frac{1}{2}s_kB_ks_k)$ 
			\STATE Obtain estimates $f_k^0$ and $f_k^s$ of $f(x_k)$ and $f(x_k + s_k)$
			\STATE Compute $\bar{g}_k = \frac{1}{b_k}\sum_{i\in O_k} \nabla f_i(x_k + s_k)$, $y_k = \bar{g}_k - g_k $
			\STATE TotalSFO = TotalSFO + $b_k$	
			\STATE Compute $\rho_k = \frac{f_k^0 - f_k^d}{Pred}$
			\IF {$\rho_k \geq \eta_1$}
			\STATE  $x_{k+1} = x_k + s_k$
			\STATE $\mu_{k+1} = \min(\gamma\mu_k, \mu_{\max})$
			\ELSE 
			\STATE $x_{k+1} = x_k$	
			\STATE $\mu_{k+1} = \mu_k/\gamma$ 
			\ENDIF			
			\IF {$\left|s_k^{T}(y_k - B_ks_k) \right| \geq r \left\|s_k \right\|\left\|y_k -B_ks_k \right\|   $}
			\STATE$S_{k+1}= \left[ S_k,s_k \right], Y_{k+1} = \left[ Y_k,y_k \right]  $	
			\IF{ len($S_{k+1}$) $\geq m $ }
			\STATE delete $S_{k+1}[1]$, $Y_{k+1}[1]$
			\ENDIF
			\ENDIF	 
			\STATE Set $k := k+1$
			\ENDWHILE
		\end{algorithmic}
	\end{algorithm}

	%\begin{acknowledgements}
	%\end{acknowledgements}
	
	\bibliographystyle{spmpsci} 
	% mathematics and physical sciences
	
	\bibliography{STRME}

	% etc

\end{document}